\documentclass[12pt]{amsart}
\usepackage{amssymb,amsmath,amsfonts,latexsym,setspace,amscd}
\usepackage{bm}
\newcommand{\bea}{\begin{eqnarray}}
\newcommand{\eea}{\end{eqnarray}}

\newcommand{\vp}{\varphi}
\newcommand{\cd}{\cdot}

\newcommand{\clb}{\mathcal{B}}

\newcommand{\cld}{\mathcal{D}}
\newcommand{\cle}{\mathcal{E}}
\newcommand{\clf}{\mathcal{F}}

\newcommand{\clh}{\mathcal{H}}
\newcommand{\clk}{\mathcal{K}}
\newcommand{\cll}{\mathcal{L}}
\newcommand{\clm}{\mathcal{M}}
\newcommand{\cln}{\mathcal{N}}
\newcommand{\clo}{\mathcal{O}}

\newcommand{\clq}{\mathcal{Q}}
\newcommand{\clr}{\mathcal{R}}
\newcommand{\cls}{\mathcal{S}}

\newcommand{\clw}{\mathcal{W}}

\newcommand{\z}{\bm{z}}
\newcommand{\w}{\bm{w}}

\newcommand{\raro}{\rightarrow}

\def \qed {\hfill \vrule height6pt width 6pt depth 0pt}
\def\textmatrix#1&#2\\#3&#4\\{\bigl({#1 \atop #3}\ {#2 \atop #4}\bigr)}
\def\dispmatrix#1&#2\\#3&#4\\{\left({#1 \atop #3}\ {#2 \atop #4}\right)}
\newcommand{\be}{\begin{equation}}
\newcommand{\ee}{\end{equation}}
\newcommand{\ben}{\begin{eqnarray*}}
\newcommand{\een}{\end{eqnarray*}}

\newcommand{\NI}{\noindent}

\newcommand{\bi}{\begin{itemize}}
\newcommand{\ei}{\end{itemize}}

\newtheorem{Theorem}{\sc Theorem}[section]
\newtheorem{Lemma}[Theorem]{\sc Lemma}
\newtheorem{Proposition}[Theorem]{\sc Proposition}
\newtheorem{Corollary}[Theorem]{\sc Corollary}
\newtheorem{Definition}[Theorem]{\sc Definition}
\newtheorem{Example}[Theorem]{\sc Example}
\newtheorem{Remark}[Theorem]{\sc Remark}
\newtheorem{Note}[Theorem]{\sc Note}
\newtheorem{Question}{\sc Question}
\newtheorem{ass}[Theorem]{\sc Assumption}
\newcommand{\bt}{\begin{Theorem}}
\def\beginlem{\begin{Lemma}}
\def\beginprop{\begin{Proposition}}
\def\begincor{\begin{Corollary}}
\def\begindef{\begin{Definition}}
\def\beginexamp{\begin{Example}}
\def\beginrem{\begin{Remark}}
\def\beginq{\begin{Question}}
\def\beginass{\begin{ass}}
\def\beginnote{\begin{Note}}
\newcommand{\et}{\end{Theorem}}
\def\endlem{\end{Lemma}}
\def\endprop{\end{Proposition}}
\def\endcor{\end{Corollary}}
\def\enddef{\end{Definition}}
\def\endexamp{\end{Example}}
\def\endrem{\end{Remark}}
\def\endq{\end{Question}}
\def\endass{\end{ass}}
\def\endnote{\end{Note}}

\begin{document}

\title[Applications of Hilbert Module Approach to Multivariable Operator Theory]{Applications of Hilbert Module Approach to Multivariable Operator Theory}

\author[Jaydeb Sarkar]{Jaydeb Sarkar}
\address{Indian Statistical Institute, Statistics and Mathematics Unit, 8th Mile, Mysore Road, Bangalore, 560059, India}
\email{jay@isibang.ac.in, jaydeb@gmail.com}




\begin{abstract}

A commuting $n$-tuple $(T_1, \ldots, T_n)$ of bounded linear
operators on a Hilbert space $\clh$ associate a Hilbert module
$\clh$ over $\mathbb{C}[z_1, \ldots, z_n]$ in the following sense:
\[\mathbb{C}[z_1, \ldots, z_n] \times \clh \raro \clh, \quad \quad (p, h) \mapsto
p(T_1, \ldots, T_n)h,\]where $p \in \mathbb{C}[z_1, \ldots, z_n]$
and $h \in \clh$. A companion survey provides an introduction to the
theory of Hilbert modules and some (Hilbert) module point of view to
multivariable operator theory. The purpose of this survey is to
emphasize algebraic and geometric aspects of Hilbert module approach
to operator theory and to survey several applications of the theory
of Hilbert modules in multivariable operator theory. The topics
which are studied include: generalized canonical models and
Cowen-Douglas class, dilations and factorization of reproducing
kernel Hilbert spaces, a class of simple submodules and quotient
modules of the Hardy modules over polydisc, commutant lifting
theorem, similarity and free Hilbert modules, left invertible
multipliers, inner resolutions, essentially normal Hilbert modules,
localizations of free resolutions and rigidity phenomenon.

\NI This article is a companion paper to ``An Introduction to
Hilbert Module Approach to Multivariable Operator Theory''.
\end{abstract}

\maketitle

\tableofcontents

\section{Introduction}\label{I}

The main motivation of Hilbert module approach to (multivariable)
operator theory is fourfold: (1) elucidating role of
Brown-Douglas-Fillmore theory (1973) to operator theory, (2) complex
geometric interpretation of (a class of) reproducing kernel Hilbert
spaces in the sense of Cowen-Douglas class (1978), (3) Hormandar's
algebraic approach, in the sense of Koszul complex, to corona
problem (1967) and (4) Taylor's notion of joint spectrum (1970),
again in the sense of Koszul complex, in operator theory and
function theory.

The general topics for this article is to survey several
applications of complex geometry and commutative algebra, with a
view of (Hilbert) module approach, to multivariable operator theory.

It is hoped that the formalism and observations presented here will
provide better understanding of the problems in operator theory in a
more general framework. The underlying idea of this survey is to:

\NI (i) Study generalized canonical models and make connections
between the multipliers and the quotient modules on one side, and
the hermitian anti-holomorphic vector bundles and curvatures on the
other side (see Section 2).

\NI (ii) Determine when a quasi-free Hilbert module can be realized
as a quotient module of a reproducing kernel Hilbert module (see
Section 3).

\NI (iii) Analyze Beurling type representation of (a class of)
submodules and quotient modules of $H^2(\mathbb{D}^n)$, $n >1$ (see
Section 4).

\NI (iv) Determine when a Hilbert module over $\mathbb{C}[\bm{z}]$
is similar to a quasi-free (or reproducing kernel) Hilbert module
(see Section 5).

\NI (v) Analyze similarity problem for generalized canonical models
corresponding to corona pairs in $H^{\infty}(\mathbb{D})$ (see
Section 6).

\NI (vi) Analyze free resolutions of Hilbert modules and
corresponding localizations and to relate with the Taylor's joint
spectrum (see Section 7).

\NI (vii) Study the rigidity properties, that is, to determine the
lattice of submodules of a reproducing kernel Hilbert module, up to
unitarily equivalence (see Section 8).

\NI (viii) Determine when a Hilbert module is small, that is, when a
(reproducing kernel) Hilbert module is essentially normal (see
Section 9).

\NI\textsf{Notations and Conventions:} (i) $\mathbb{N} =$ Set of all
natural numbers including 0.  (ii) $n \in \mathbb{N}$ and $n \geq
1$, unless specifically stated otherwise. (iii) $\mathbb{N}^n =
\{\bm{k} = (k_1, \ldots, k_n) : k_i \in \mathbb{N}, i = 1, \ldots,
n\}$. (iv)  $\mathbb{C}^n =$ the complex $n$-space. (v) $\Omega$ :
Bounded domain in $\mathbb{C}^n$. (vi) $\bm{z} = (z_1, \ldots, z_n)
\in \mathbb{C}^n$. (vii) ${z}^{\bm{k}} = z_1^{k_1}\cdots z_n^{k_n}$.
(viii) $\clh, \clk, \cle, \cle_*$ : Hilbert spaces. (ix) $\clb(\clh,
\clk)=$ the set of all bounded linear operators from $\clh$ to
$\clk$. (x) $T = (T_1, \ldots, T_n)$, $n$-tuple of commuting
operators. (xi) $T^{\bm{k}} = T_1^{k_1} \cdots T_n^{k_n}$. (xii)
$\mathbb{C}[\z] = \mathbb{C}[z_1, \ldots, z_n]$. (xiii)
$\mathbb{D}^n = \{\z : |z_i| <1, i = 1, \ldots, n\}$, $\mathbb{B}^n
= \{\z : \|\z\|_{\mathbb{C}^n} <1\}$. (xiv) $H^2_{\cle}(\mathbb{D})$
: $\cle$-valued Hardy space over $\mathbb{D}$.

\NI Throughout this note all Hilbert spaces are over the complex
field and separable. Also for a closed subspace $\cls$ of a Hilbert
space $\clh$, the orthogonal projection of $\clh$ onto $\cls$ will
be denoted by $P_{\cls}$.

This article is a companion paper to ``An Introduction to Hilbert
Module Approach to Multivariable Operator Theory'' (see
\cite{JS-HB}).

\section{Generalized canonical models in the Cowen-Douglas
class}\label{GCMCDC}

Let $\cle$ and $\cle_*$ be Hilbert spaces and $\clh \in
B_1^*(\Omega)$. Moreover, assume $\Theta \in \clm_{\clb(\cle,
\cle_*)}(\clh)$. Then the quotient module $\clh_{\Theta} = \clh
\otimes \cle_*/ \Theta(\clh \otimes \cle)$ is called the
\textit{generalized canonical model} associated with $\clh$ and
$\Theta$. In other words, a generalized canonical model can be
obtained by the resolution

\[ \cdots \longrightarrow \clh
\otimes \cle \stackrel{M_{\Theta}} \longrightarrow \clh \otimes
\cle_* \stackrel{\pi_{\Theta}} \longrightarrow \clh_{\Theta}
\longrightarrow 0.\]This is a generalization of Sz.-Nagy-Foias
notion of canonical model (see Section 4 in \cite{JS-HB}) to
quotient modules of Hilbert modules.

Let $\clh \in B_1^*(\mathbb{D})$  be a contractive Hilbert module
over $A(\mathbb{D})$. Then $\clh$ is in $C_{\cdot 0}$ class and the
characteristic function $\Theta_{\clh}$, in the sense of Sz.-Nagy
and Foias, is a complete unitary invariant (see Section 4 in
\cite{JS-HB}). On the other hand, the curvature, in the sense of
Cowen and Douglas, is another complete unitary invariant. A very
natural question then arises: whether the characteristic function is
connected with the curvature of the canonical model of $\clh$.

One can formulate the above problem in a more general framework by
replacing the Hardy module with a Hilbert module in $B_1^*(\Omega)$.
More precisely, let $\clh \in B^*_1(\Omega)$ and $\Theta \in
\clm_{\clb(\cle, \cle_*)}(\clh)$. Suppose the quotient module
$\clh_{\Theta} = \clh \otimes \cle_*/ \Theta(\clh \otimes \cle)$ is
in $B_m^*(\Omega)$. Does there exists any connection between the
multipliers and curvature corresponding to the hermitian
anti-holomorphic vector bundle $E_{\clh_{\Theta}}$?

The purpose of this section is to study generalized canonical models
and make connections between the multipliers and the quotient
modules on one side, and the hermitian anti-holomorphic vector
bundles and curvatures on the other side. Results concerning
similarity and unitarily equivalence will be derived from these
connections. The final subsection of this section will discuss some
quotient modules of the familiar Hardy and weighted Bergman modules
over $A(\mathbb{D})$ and trace basic facts about unitary equivalence
and curvature equality.

\subsection{Generalized canonical models in $B^*_m(\Omega)$}

Generalized canonical models yields a deeper understanding of many
issues in the study of Hilbert modules. However, the present
approach we will assume only finite dimensional coefficient spaces
with left invertible multiplier:
\[ 0 \longrightarrow \clh
\otimes \mathbb{C}^p \stackrel{M_{\Theta}} \longrightarrow \clh
\otimes \mathbb{C}^{q} \stackrel{\pi_{\Theta}} \longrightarrow
\clh_{\Theta} \longrightarrow 0,\]where $p, q \in \mathbb{N}$ and $q
> p$.
\begin{Theorem}\label{n-curv}
Let $1 \leq p < q$ and $\clh_{\Theta}$ be a generalized canonical
model corresponding to $\clh \in B^*_1(\Omega)$ and a left
invertible $\Theta \in \clm_{\clb(\mathbb{C}^p,
\mathbb{C}^q)}(\clh)$. Then

(1) $\clh_{\Theta} \in B_{q-p}^*(\Omega)$, and

(2) $V^*_\Theta(\w)= (\mbox{ran}\text{ }\Theta(\w))^\perp =
\mbox{ker} \text{ }\Theta(\w)^*$ defines a hermitian
anti-holomorphic vector bundle \[V_{\Theta}^* =\coprod_{\w \in
\Omega} V_\Theta(\w)^*,\] over $\Omega$ such that
\[E^*_{\clh_{\Theta}} \cong E^*_{\clh} \otimes V^*_{\Theta}.\]

\NI In particular, if $q = p+1$ then $\clh_{\Theta} \in
B_1^*(\Omega)$ and $V_{\Theta}$ is a line bundle.
\end{Theorem}

\NI\textsf{Proof.}  Localizing the short exact sequence of Hilbert
modules
\[0 \longrightarrow \clh  \otimes \mathbb{C}^p
\stackrel{M_{\Theta}}\longrightarrow \clh \otimes \mathbb{C}^{q}
\stackrel{\pi_{\Theta}}\longrightarrow \clh_{\Theta} \longrightarrow
0,\] at $\w \in \Omega$, that is, taking quotients by $I_{\w} \cdot
(\clh \otimes \mathbb{C}^p)$,  $I_{\w} \cdot (\clh \otimes
\mathbb{C}^{q}),$ and $I_{\w} \cdot \clh_{\Theta}$, respectively,
one obtain the following exact sequence (see Theorem 5.12 in
\cite{DP})\[\mathbb{C}_{\w} \otimes \mathbb{C}^p
\stackrel{I_{\mathbb{C}_{\w}} \otimes \Theta(\w)} \longrightarrow
\mathbb{C}_{\w} \otimes \mathbb{C}^{q} \stackrel{\pi_{\Theta}(\w)}
\longrightarrow \clh_{\Theta}/ I_{\w} \cdot \clh_{\Theta}
\longrightarrow 0.\] Since $\mbox{dim}[\mbox{ran } \Theta(\w)] = p$
for all $\w \in \Omega$, it follows that $\mbox{dim} \Big[\mbox{ker
} \pi_{\Theta}(\w)\Big] = p$, and thus
\[\mbox{dim~} \Big[\clh_{\Theta}/ I_{\w} \cdot \clh_{\Theta}\Big] = \mbox{dim~}
\Big[\clh_{\Theta}\Big/\Big(\sum_{i=1}^n (M_{z_i} - w_i I_{\clh})
\clh_{\Theta}\Big)\Big] = q-p,\] that is, \[\mbox{dim}
\Big[\mathop{\cap}_{i=1}^n \mbox{ker~} (M_{z_i} - w_i
I_{\clh})^*|_{\clh_{\Theta}}\Big] = q-p,\]for all $\w \in \Omega$.

\NI The next step is to prove the following equality
\[\bigvee_{w \in \Omega} \{\mbox{ker~}(M_z - w I_{\clh})^* \otimes
{\mbox{ker~}\Theta(w)^*}\} = (\clh \otimes \mathbb{C}^q) \ominus
\mbox{ran~}M_{\Theta}.\] For simplicity of notation, assume that $q
= p+1$. The proof of the general case is essentially the same as the
one presented below (or see Theorem 3.3 in \cite{DKKS2}). To this
end, let $\{e_i\}_{i=1}^{p+1}$ be the standard orthonormal basis for
$\mathbb{C}^{p+1}$ and let $\Delta_{\Theta}$ be the formal
determinant
\[
\Delta_{\Theta}(\w) = \mbox{det}\
\begin{bmatrix} e_1 & \theta_{1,1}(\w) & \cdots & \theta_{1,p}(\w) \\ \vdots & \vdots & \vdots & \vdots
\\ e_{p+1} & \theta_{p+1,1}(\w) & \cdots & \theta_{p+1,p}(\w) \end{bmatrix} \in \mathbb{C}^{p+1},
\]
where $\Theta(\w) = (\theta_{i,j}(\w))$ and $\w \in \Omega$. Since
$\Theta(\w)$ has a left inverse $\Psi(\w)$, it follows that
$\mbox{rank~} \Theta(\w)= l$, and hence $\Delta_{\Theta}(\w) \neq 0$
for all $\w \in \Omega$. Set $\gamma_{\w} := k_{\w} \otimes
\overline{\Delta_{\Theta}(\w)} \neq 0$ for all $\w \in \Omega$,
where $k_{\w}$ is any non-zero vector in $E^*_{\clh}(\w) \subseteq
\clh$ and $\overline{\Delta_{\Theta}(\w)}$ is the complex conjugate
of $\Delta_{\Theta}(\w)$ relative to the basis
$\{e_i\}_{i=1}^{p+1}$. Moreover, consider the inner product of
$\gamma_{\w}$ with
\[M_{\Theta} \begin{bmatrix}h_1\\ \vdots\\h_l\end{bmatrix} =
\begin{bmatrix}\sum_{j=1}^p \theta_{1,j} h_j\\\vdots\\ {\sum}_{j=1}^p
\theta_{p+1, j} h_j\end{bmatrix} \in \clh \otimes \mathbb{C}^{p+1},
\]for $\{h_i\}_{i=1}^p \subseteq \clh$. Evaluating the resulting
functions at $\w \in \Omega$, one can conclude that these functions
are the sum of the products of $h_i(\w)$ with coefficients equal to
the determinants of matrices with repeated columns and hence
\[\langle M_{\Theta} \begin{bmatrix}h_1\\ \vdots\\h_p\end{bmatrix},
\gamma_{\w} \rangle = 0. \] Thus, $\gamma_{\w} \perp
\mbox{ran~}M_\Theta$ for all $\w \in \Omega$. Also, it is easy to
see that
\[(M_{z_i}^* \otimes I_{\mathbb{C}^{p+1}}) \gamma_{\w} = \bar{w_i}
\gamma_{\w}, \]for $\w \in \Omega$ and for all $i = 1, \ldots, n$,
so that
\[\mathop{\bigcap}_{i=1}^n \mbox{ker~} (M_{z_i} \otimes I_{\mathbb{C}^{p+1}} - w_i I_{\clh \otimes \mathbb{C}^{p+1}})^*|_{\clh_{\Theta}} = \mathbb{C}\cdot
\gamma_{\w},\]for all $\w \in \Omega$.

\NI The next step is to prove that $\bigvee_{\w \in \Omega}k_{\w}
\otimes \overline{\Delta_{\Theta}(w)}= {\clh_{\Theta}}$. For all $g
=\sum_{i=1}^{p+1} g_i \otimes e_i \in \clh \otimes \mathbb{C}^{p+1}$
with $g \perp \gamma_{\w}$ for every $\w \in \Omega$, one must
exhibit the representation $g_i(\w) = \sum_{j=1} ^p \eta_j(\w)
\theta_{ij}(\w)$ for $i=1,...,p+1$, where the $\{\eta_j\}_{j=1}^p$
are functions in $\clh$. Fix $\w_0 \in \Omega$. The assumption
$\langle g, \gamma_{\w_0} \rangle =0$ implies that
\begin{equation}\label{3.1}
\mbox{det}\ \begin{bmatrix} g_1(\w_0) & \theta_{1,1}(\w_0) & \cdots
& \theta_{1,p}(\w_0) \\ \vdots & \vdots & \vdots & \vdots \\
g_{p+1}(\w_0) & \theta_{p+1,1}(\w_0) & \cdots & \theta_{p+1,p}(\w_0)
\end{bmatrix}=0.
\end{equation}
Now view the matrix
$$\Theta(\w_0)=\begin{bmatrix} \theta_{1,1}(\w_0) & \cdots &
\theta_{1,p}(\w_0) \\ \vdots & \vdots & \vdots \\
\theta_{p+1,1}(\w_0) & \cdots & \theta_{p+1,l}(\w_0) \end{bmatrix}$$
as the coefficient matrix of a linear system of $(p+1)$ equations in
$p$ unknowns. Since $\mbox{rank~}\Theta(\w_0)= p$, some principal
minor (which means taking some $p$ rows) has a non-zero determinant.
Hence, using Cramer's rule, one can uniquely solve for
$\{\eta_j(\w_0)\}_{j=1}^p \subseteq \mathbb{C}^p$, at least for
these $p$ rows. But by (\ref{3.1}), the solution must also satisfy
the remaining equation. Hence we obtain the
$\{\eta_j(\w_0)\}_{j=1}^p \subseteq \mathbb{C}^p$ and define
$$\xi(\w_0)=\sum_{j=1}^{p} \eta_j(\w_0) \otimes e_j,$$
so that
$$
g(\w_0)=\Theta(\w_0)\xi(\w_0),
$$
for each $\w_0 \in \Omega$. After doing this for each $\w \in
\Omega$, we use the left inverse $\Psi(\w)$ for $\Theta(\w)$ to
obtain
$$\xi(\w) =(\Psi(\w) \Theta(\w))\xi(\w)=\Psi(\w)(\Theta(\w)
\xi(\w))=\Psi(\w) g(\w) \in \clh \otimes \mathbb{C}^p.$$
Consequently, $\{\eta_j\}_{j=1}^p \subseteq \clh$ and $\bigvee_{\w
\in \Omega} \gamma_{\w} = \clh_{\Theta}$.

\NI Lastly, the closed range property of $\clh_{\Theta}$ follows
from that of $\clh$. In particular, since the column operator $M_z^*
- \bar{w}I_{\clh}$ (see Definition 3.1 in \cite{JS-HB}) acting on
$\clh \otimes \mathbb{C}^{l+1}$ has closed range and a finite
dimensional kernel, it follows that restricting it to the invariant
subspace $\clh_{\Theta} \subseteq \clh \otimes \mathbb{C}^{p+1}$
yields an operator with closed range and hence $\clh_{\Theta} \in
B_1^*(\Omega)$. \qed

The above result allows one to construct a wide range of
Cowen-Douglas Hilbert modules over domains in $\mathbb{C}^n$.

\subsection{Curvature equality}
The following is a very useful equality for the class of generalized
canonical models.
\begin{Theorem}\label{n-curv-equal}

Let $1 \leq p < q$ and $\clh_{\Theta}$ be a generalized canonical
model corresponding to $\clh \in B^*_1(\Omega)$ and a left
invertible $\Theta \in \clm_{\clb(\mathbb{C}^p,
\mathbb{C}^q)}(\clh)$. Then \[\clk_{E^*_{\clh_{\Theta}}} -
\clk_{E^*_{\clh}} = \clk_{V^*_{\Theta}}.\]
\end{Theorem}
\NI\textsf{Proof.} To establish the curvature formula, first recall
that the formula for the curvature of the Chern connection on an
open subset $U \subseteq \Omega$ for a hermitian anti-holomorphic
vector bundle is $\bar{\partial}[G^{-1}{\partial}G]$, where $G$ is
the Gramian for an anti-holomorphic frame $\{f_i\}_{i=1}^{q-p}$ for
the vector bundle on $U$ (cf. \cite{CS}). Assume that $U$ is chosen
so that the $\{k_{\w}\}$ for $\w \in \Omega$ can be chosen to be an
anti-holomorphic function on $U$. Denoting by $G_{\Theta}$ the
Gramian for the frame $\{k_{\w} \otimes f_i(\w)\}_{i=1}^{q-p}$,
$G_{\Theta}(\w)$ equals the $(q-p) \times (q-p)$ matrix
\[G_{\Theta}(\w) = \big( \langle k_{\w} \otimes f_i(\w), k_{\w} \otimes
f_j(\w)\rangle \big)_{i, j = 1}^{q-p} = \|k_{\w}\|^2 \big ( \langle
f_i(\w), f_j(\w)\rangle \big)_{i,j = 1}^{q-p} = \|k_{\w}\|^2
G_{f}(\w),\]where $G_{f}$ is the Gramian for the anti-holomorphic
frame $\{f_i(\w)\}_{i=1}^{q-p}$ for $V^*_\Theta$. Then
\[
\begin{split}
\bar{\partial} [G_{\Theta}^{-1}({\partial} G_{\Theta})] &=
\bar{\partial} [\frac{1}{\|k_{\w}\|^2}
G_{f}^{-1}({\partial}(\|k_{\w}\|^2 G_{f}))]
\\& =
\bar{\partial} [\frac{1}{\|k_{\w}\|^2}
G_{f}^{-1}({\partial}(\|k_{\w}\|^2) G_{f} + \|k_{\w}\|^2 {\partial}
G_{f})] \\& = \bar{\partial}
[\frac{1}{\|k_{\w}\|^2}{\partial}(\|k_{\w}\|^2) + G_{f}^{-1}
{\partial} G_{f} ]\\ & = \bar{\partial}
[\frac{1}{\|k_{\w}\|^2}{\partial}(\|k_{\w}\|^2)] + \bar{\partial} [
G_{f}^{-1}{\partial} G_{f}].
\end{split}
\]
Hence, expressing these matrices in terms of the respective frames
and using the fact that the coordinates of a bundle and of its dual
can be identified using the basis given by the frame, one has
\[\clk_{E^*_{\clh_{\Theta}}}(\w) - \clk_{E^*_{\clh}}(\w) \otimes
I_{V^*_\Theta(\w)} = I_{E^*_\clh (\w)} \otimes
\clk_{V^*_\Theta}(\w),\] for all $\w \in U$. Since the coordinate
free formula does not involve $U$, this completes the proof. \qed

Based on Theorems \ref{n-curv} and \ref{n-curv-equal}, one can say
that the isomorphism of quotient Hilbert modules is independent of
the choice of the basic Hilbert module "building blocks" from which
they were created.

\begin{Corollary}\label{n-curv-cor}
Let $\clh, \tilde{\clh} \in B_1^*(\Omega)$ and $\Theta_1, \Theta_2
\in \clm_{\clb(\mathbb{C}^p, \mathbb{C}^q)} (\clh) \cap
\clm_{\clb(\mathbb{C}^p, \mathbb{C}^q)} (\tilde{\clh})$ are left
invertible with inverse in $\clm_{\clb(\mathbb{C}^q, \mathbb{C}^p)}
(\clh) \cap \clm_{\clb(\mathbb{C}^q, \mathbb{C}^p)} (\tilde{\clh})$.
Then $\clh_{\Theta_1}$ is isomorphic to $\clh_{\Theta_2}$ if and
only if $\tilde{\clh}_{\Theta_1}$ is isomorphic to
$\tilde{\clh}_{\Theta_2}$.
\end{Corollary}

\NI\textsf{Proof.} The statement is obvious from the tensor product
representations $E^*_{\clh_{\Theta_i}} \cong E^*_{\clh} \otimes
V^*_{\Theta_i}$ and $E^*_{\tilde{\clh}_{\Theta_i}} \cong
E^*_{\tilde{\clh}} \otimes V^*_{\Theta_i}$, for $i=1, 2$ ; that is,
isomorphic as hermitian anti-holomorphic bundles, and the result
that $\clk_{E^*_{\clh_{\Theta_1}}} = \clk_{E^*_{\clh_{\Theta_2}}}$
if and only if $\clk_{V^*_{\Theta_1}} = \clk_{V^*_{\Theta_2}}$ as
two forms. \qed

In what follows, $\bigtriangledown^2$ denotes the Laplacian
\[\bigtriangledown^2 = 4 \partial
\bar{\partial} = 4 \bar{\partial}{\partial}.\]

\begin{Theorem}\label{cor-3.2}
Let $\clh \in B_1^*(\Omega)$ and $\Theta_1, \Theta_2 \in
\clm_{\clb(\mathbb{C}^l, \mathbb{C}^{l+1})}(\clh)$ are left
invertible multipliers. Then the quotient Hilbert modules
$\clh_{\Theta_1}$ and $\clh_{\Theta_2}$ are isomorphic if and only
if \[ \bigtriangledown^2 \mbox{log~} \|\Delta_{\Theta_1}\| =
\bigtriangledown^2 \mbox{log~}\|\Delta_{\Theta_2}\|,\]where
$\Delta_{\Theta_i}$ is an anti-holomorphic cross section of
$V^*_{\Theta_i}$ and $i = 1, 2$.
\end{Theorem}

\NI\textsf{Proof.} Choose a cross section $k_{\w}$ so that $k_{\w}
\otimes \overline{\Delta_{\Theta_i}(\w)}$, $i=1, 2$, are
anti-holomorphic local cross-sections of $E^*_{\clh_{\Theta_1}}$ and
$E^*_{\clh_{\Theta_2}}$, respectively, over some open subset $U
\subseteq \Omega$. Since every $\w_0 \in \Omega$ is contained in
such an open subset $U$ of $\Omega$, by rigidity theorem \cite{CD}
(or Theorem 3.2 in \cite{JS-HB}), it follows that $\clh_{\Theta_1}
\cong \clh_{\Theta_2}$ if and only if
$$\clk_{E^*_{\clh_{\Theta_1}}}(\z) = \clk_{E^*_{\clh_{\Theta_2}}}(\z),$$
for every $\z \in \Omega$ or, equivalently, \[\bigtriangledown^2
\mbox{log~} \|\Delta_{\Theta_1}\| = \bigtriangledown^2
\mbox{log~}\|\Delta_{\Theta_2}\|,\]by Theorem  \ref{n-curv-equal}.
This completes the proof. \qed

\subsection{Examples and applications}

The purpose of this subsection is to describe a class of simple
examples of generalized canonical models in $B^*_1(\mathbb{D})$.

Let $\Theta \in H ^{\infty} _{\clb(\mathbb{C},
\mathbb{C}^2)}(\mathbb{D})$ so that \[\Theta(z)=
\begin{bmatrix}\theta_1(z)\\ \theta_2(z)\end{bmatrix},\]and $\theta_1, \theta_2 \in
H^\infty(\mathbb{D})$ and $z \in \mathbb{D}$. $\Theta$ is said to
satisfy the {\it corona condition} if there exists an $\epsilon > 0$
such that $|\theta_1(z)|^2 + |\theta_2(z)|^2 > \epsilon$ for all $z
\in \mathbb{D}$ (see Section \ref{SFHM}).

For the rest of this subsection, fix a corona pair $\Theta =
\begin{bmatrix}\theta_1\\ \theta_2\end{bmatrix}\in H ^{\infty}
_{\clb(\mathbb{C}, \mathbb{C}^2)}(\mathbb{D})$ and use the notation
$\clh$ to denote the Hardy, the Bergman, or a weighted Bergman
module over $\mathbb{D}$. Consider the generalized canonical model
$\clh_\Theta$ corresponding to the exact sequence of Hilbert
modules:
\[0 \longrightarrow \clh  \otimes \mathbb{C}
\stackrel{M_{\Theta}}\longrightarrow \clh \otimes \mathbb{C}^2
\stackrel{\pi_{\Theta}}\longrightarrow \clh_{\Theta} \longrightarrow
0,\] where the first map $M_{\Theta}$ is $M_{\Theta} f = \theta_1 f
\otimes e_1 +\theta_2 f \otimes e_2$ and the second map
$\pi_{\Theta}$ is the quotient Hilbert module map.

Note that by taking the kernel functions for $H^2(\mathbb{D})$ and
$L^2_{a,\alpha}(\mathbb{D})$ as an anti-holomorphic cross section of
bundles an easy computation shows that
$$\clk_{E^*_{H^2(\mathbb{D})}}(z) = - \frac{1}{(1 - |z|^2)^{2}},$$ and
$$\clk_{E^*_{L^2_{a,\alpha}(\mathbb{D})}}(z) = - \frac{2+\alpha}{(1 - |z|^2)^{2}}.
$$

The following is immediate consequence of Theorems \ref{n-curv} and
\ref{n-curv-equal}.

\begin{Theorem}\label{hardy-section}
For $\Theta=\begin{bmatrix}\theta_1\\ \theta_2\end{bmatrix}$
satisfying the corona condition, $\clh_{\Theta} \in
B_1^*(\mathbb{D})$ and
\begin{equation}\label{4.2} \clk_{E^*_{{\clh}_{\Theta}}}(w)=\clk
_{E^*_{\clh}}(w)-\frac{1}{4} \bigtriangledown^2
\mbox{log}\,(|\theta_1(w)|^2 + |\theta_2(w)|^2). \quad \quad (w \in
\mathbb{D})
\end{equation}
\end{Theorem}

\begin{Theorem}
Let $\Theta=\begin{bmatrix}\theta_1\\ \theta_2\end{bmatrix}$ and
$\Phi=\begin{bmatrix}\vp_1\\ \vp_2\end{bmatrix}$ satisfy the corona
condition. The quotient Hilbert modules ${\clh}_{\Theta}$ and
${\clh}_{\Phi}$ are isomorphic if and only if
$$
\bigtriangledown^2 \mbox{log}\,\frac{|\theta_1(z)|^2 +
|\theta_2(z)|^2}{|\vp_1(z)|^2 + |\vp_2(z)|^2}=0. \quad \quad (z \in
\mathbb{D})
$$
\end{Theorem}

\NI \textsf{Proof.} Since ${\clh}_{\Theta}, {\clh}_{\Phi} \in
B_1^*(\mathbb{D})$, they are isomorphic if and only if
$\clk_{E^*_{{\clh}_{\Theta}}}(w)=\clk_{E^*_{{\clh}_{\Phi}}}(w)$ for
all $w \in \mathbb{D}$. But note that (\ref{4.2}) and an analogous
identity for $\Phi$ hold, where the $\theta_i$ are replaced with the
$\vp_i$. Since both $\Theta$ and $\Phi$ satisfy the corona
condition, the result then follows.
\qed \\

\begin{Theorem}
Suppose that $\Theta=\begin{bmatrix}\theta_1\\
\theta_2\end{bmatrix}$ and $\Phi=\begin{bmatrix}\vp_1\\
\vp_2\end{bmatrix}$ satisfy the corona condition. The quotient
Hilbert modules $(L^2_{a,\alpha}(\mathbb{D}))_{\Theta}$ and
$(L^2_{a, \beta}(\mathbb{D}))_{\Phi}$ are isomorphic if and only if
$\alpha=\beta$ and
\[
\bigtriangledown^2 \mbox{log}\,\frac{|\theta_1(z)|^2 +
|\theta_2(z)|^2}{|\vp_1(z)|^2 + |\vp_2(z)|^2}=0. \quad \quad (z \in
\mathbb{D})
\]
\end{Theorem}

\NI \textsf{Proof.} Since
$$\clk_{E^*_{(L^2_{a,\alpha}(\mathbb{D}))_\Theta}}(w) =  - \frac{2+\alpha}{(1 -
|w|^2)^{2}} - \frac{1}{4} \bigtriangledown^2
\mbox{log}\,(|\theta_1(w)|^2 + |\theta_2(w)|^2),$$ and
$$\clk_{E^*_{(L^2_{a, \beta}(\mathbb{D}))_\Phi}} (w) =  - \frac{2+\beta}{(1 -
|w|^2)^{2}} - \frac{1}{4} \bigtriangledown^2
\mbox{log}\,(|\vp_1(w)|^2 + |\vp_2(w)|^2),$$ by (\ref{4.2}), one
implication is obvious. For the other one, suppose that $(L^2_{a,
\alpha}(\mathbb{D}))_{\Theta}$ is isomorphic to $(L^2_{a,
\beta}(\mathbb{D}))_\Phi$ so that the curvatures coincide. Observe
next that $$\frac{4(\beta - \alpha)}{(1 - |w|^2)^{2}} =
\bigtriangledown^2 \mbox{log}\,\frac{|\theta_1(w)|^2 +
|\theta_2(w)|^2}{|\vp_1(w)|^2 + |\vp_2(w)|^2}.$$ Since a function
$f$ with $\bigtriangledown^2 f(z) = \frac{1}{(1 - |z|^2)^2}$ for all
$z \in \mathbb{D}$ is necessarily unbounded, one arrives at a
contradiction, unless $\alpha=\beta$ (see Lemma \ref{4.3} below).
This is due to the assumption that the bounded functions $\Theta$
and $\Phi$ satisfy the corona condition. \qed

\begin{Lemma}\label{4.3}
There is no bounded function $f$ defined on the unit disk
$\mathbb{D}$ that satisfies $\bigtriangledown^2 f(z) =\frac{1}{(1 -
|z|^2)^2}$ for all $z \in \mathbb{D}$.
\end{Lemma}

\NI\textsf{Proof.} Suppose that such $f$ exists. Since $\frac{1}{4}
\bigtriangledown^2 [(|z|^2)^{m}] = \partial \bar{\partial}
[(|z|^2)^{m}] = m^2 (|z|^2)^{m-1}$ for all $m \in \mathbb{N}$, one
see that for $$g(z) := \frac{1}{4} \sum_{m=1}^{\infty}
\frac{|z|^{2m}}{m} = - \frac{1}{4} \mbox{log}\, (1 - |z|^2),$$ $
\bigtriangledown^2 g(z) = \frac{1}{(1-|z|^2)^2}$ for all $z \in
\mathbb{D}$. Consequently, $f(z)=g(z)+h(z)$ for some harmonic
function $h$. Since the assumption is that $f$ is bounded, there
exists an $M>0$ such that $|g(z) + h(z)| \leq M$ for all $z \in
\mathbb{D}$. It follows that $$ \mbox{exp}\,(h(z))  \leq
\mbox{exp}\,( - g(z) + M) = (1 - |z|^2)^{\frac{1}{4}}
\,\mbox{exp}\,(M),$$ and letting $z = r e^{i\theta},$ we have
$\mbox{exp}\,(h(r e^{i\theta})) \leq (1 - r^2)^{\frac{1}{4}}
\mbox{exp}\,(M)$. Thus $\mbox{exp}\,(h(r e^{i\theta})) \raro 0$
uniformly as $r \raro 1^-$, and hence $\exp h(z) \equiv 0$. This is
due to the maximum modulus principle because $\exp h(z) = |
\exp(h(z) + i \tilde{h}(z))|$, where $\tilde{h}$ is a harmonic
conjugate for $h$. This leads to a contradiction, and the proof is
complete. \qed

\begin{Theorem}
For $\Theta=\begin{bmatrix}\theta_1\\
\theta_2\end{bmatrix}$ and $\Phi=\begin{bmatrix}\vp_1\\
\vp_2\end{bmatrix}$ satisfying the corona condition,
$(H^2(\mathbb{D}))_{\Theta}$ cannot be isomorphic to $(L^2_{a,
\alpha}(\mathbb{D}))_\Phi$.
\end{Theorem}

\NI \textsf{Proof.} By identity (4.2), one can conclude that
$(H^2)_{\Theta}$ is isomorphic to $(A^2_{\alpha})_{\Phi}$ if and
only if $$\frac{4(1+\alpha)}{(1 - |w|^2)^{2}} =  \bigtriangledown^2
\mbox{log}\,\frac{|\vp_1(w)|^2 + |\vp_2(w)|^2}{|\theta_1(w)|^2 +
|\theta_2(w)|^2}.$$ But according to Lemma 4.6, this is impossible
unless $\alpha=-1$.\qed

\NI \textbf{Further results and comments:}

\begin{enumerate}

\item
Let $\cle$ and $\cle_*$ be two Hilbert spaces and $\Theta \in \clo(
\Omega, \clb(\cle, \cle_*))$. One can define the holomorphic kernel
and co-kernel bundles with fibers $\mbox{ker}\, \Theta(\w)$ and
$\mbox{coker}\, \Theta(\w) = \cle_*/ \Theta(\w) \cle$ for $\w \in
\Omega$, respectively, whenever it make sense. Moreover, related
Hilbert modules with $\clh \in B^*_m(\Omega)$ can be defined for an
arbitrary $m \geq 1$. Here consideration is restricted to the
``simplest'' case, when $\Theta$ is left invertible, and obtain some
of the most ``direct'' possible results.

\item Let $\clh \in B_m^*(\mathbb{D})$ be a contractive Hilbert
module over $A(\mathbb{D})$. Then one can prove that $\clh$ is in
the $C_{\cdot 0}$ class. In this case, the connection between the
characteristic function $\Theta_{\clh}$ and the curvature of the
generalized canonical model, that is, the Sz.-Nagy-Foias canonical
model $H^2_{\cld_*}/ \Theta_{\clh} H^2_{\cld}(\mathbb{D})$, was
addressed earlier by Uchiyama in \cite{U}. His theory is
instrumental in the study of generalized canonical models (cf.
\cite{KT}, \cite{S-CD}).

\item All results presented in this section can be found in
\cite{DKKS1} and \cite{DKKS2}.

\item In connection with this section, see also the work by Zhu
\cite{Z-00}, Eschmeier and Schmitt \cite{ES} and Kwon and Treil
\cite{KT} and Uchiyama \cite{U} (see also \cite{S-CD}).

\end{enumerate}

\section{Dilation to quasi-free Hilbert modules}\label{DQFHM}

Recall that a Hilbert module $\clh$ over $\mathbb{C}[z]$ is
$C_{\cdot 0}$-contractive if and only if (see Section 4 in
\cite{JS-HB}) there exists a resolution of Hilbert modules \[0
\longrightarrow \clf_1 \stackrel{i} \longrightarrow \clf_2
\stackrel{\pi} \longrightarrow \clh \longrightarrow 0,\]where
$\clf_i = H^2_{\cle_i}(\mathbb{D})$ for some Hilbert spaces $\cle_1$
and $\cle_2$.

Now let $\clh$ be a $C_{\cdot 0}$-contractive Hilbert module over
$\mathbb{C}[\z]$  (that is, $M_i \in C_{\dot 0}$ for each $i$) and
$n \geq 2$. If one attempts to obtain a similar resolution for
$\clh$, then one quickly runs into trouble. In particular, if $n >
2$ then Parrott's example \cite{P} shows that, in general, an
isometric dilation need not exist. On the other hand, a pair of
commuting contractions is known to have an isometric dilation
\cite{An}, that is, a resolution exists for contractive Hilbert
module over $\mathbb{C}[z_1, z_2]$. However, such dilations are not
necessarily unique, that is, one can not expect that $\clf_2$ to be
a free module $H^2(\mathbb{D}^2) \otimes \cle_2$.

The purpose of this section is to study the following problem: Let
$\clr \subseteq \clo(\Omega, \mathbb{C})$ be a reproducing kernel
Hilbert module over $A(\Omega)$ and $\clm$ be a quasi-free Hilbert
module over $A(\Omega)$. Determine when $\clm$ can be realized as a
quotient module of the free module $\clr \otimes \cle$ for some
coefficient space $\cle$, that is, when $\clm$ admits a free
resolution \[0 \longrightarrow \cls \stackrel{i} \longrightarrow
\clr \otimes \cle \stackrel{\pi} \longrightarrow \clm
\longrightarrow 0,\]where $\cls$ is a submodule of $\clr \otimes
\cle$.

Another important motivation for studying dilation to quasi-free
Hilbert modules is to develop some connections between free
resolutions, positivity of kernel functions and factorizations of
kernel functions. Our main tool is to establish a close relationship
between the kernel functions for the Hilbert modules in an exact
sequence using localization.

\subsection{Factorization of reproducing kernels}

Let $\clr \subseteq \clo(\Omega, \mathbb{C})$ be a reproducing
kernel Hilbert space and $\clh$ be a quasi-free Hilbert module of
multiplicity $m$ over $\mathbb{C}[\bm{z}]$ or $A(\Omega)$ and $\cle$
a Hilbert space. Then $\clr \otimes \cle$ being a dilation of $\clh$
is equivalent to the exactness of the sequence of Hilbert modules
$$0 \longrightarrow  \cls \stackrel{i} \longrightarrow \clr \otimes \cle
\stackrel{\pi} \longrightarrow \clh \longrightarrow 0,$$where the
second map is the inclusion $i$ and the third map is the quotient
map $\pi$ which is a co-isometry. The aim of this subsection is to
relate the existence of an $\clr \otimes \cle$-dilation of a
reproducing kernel Hilbert module $\clh_K$ to the positivity of the
kernel function $K$.

\begin{Theorem}\label{TH1}
Let $\clr \subseteq \clo(\Omega, \mathbb{C})$ be a reproducing
kernel Hilbert module with the scalar kernel function $k$ and $\clh$
be a quasi-free Hilbert module of multiplicity $m$  over $A(\Omega)$
or $\mathbb{C}[\bm{z}]$. Then $\clr \otimes \cle$ is a dilation of
$\clh$ for some Hilbert space $\cle$, if and only if there is a
holomorphic map $\pi_{\bm{z}} \in \clo(\Omega, \cll(\cle, l^2_m))$
such that
$$K_{\clh}(\bm{z}, \bm{w}) = k(\bm{z}, \bm{w}) \pi_{\bm{z}}
\pi_{\bm{w}}^*. \quad \quad (\bm{z}, \bm{w} \in \Omega)$$
\end{Theorem}

\NI \textsf{Proof.} Let $\clr \otimes \cle$ be a dilation of $\clh$,
that is,\[0 \raro \cls \raro \clr \otimes \cle \raro \clh \raro 0.\]
Localizing the above exact sequence of Hilbert modules at $\bm{z}
\in \Omega$ one arrives at \setlength{\unitlength}{3mm}
 \begin{center}
 \begin{picture}(40,16)(0,0)
\put(2,3){$ \cls/I_{\bm{z}} \cls$} \put(10,3){$ (\clr \otimes
\cle)/I_{\bm{z}} (\clr \otimes \cle) $} \put(26,3){$ \clh/I_{\bm{z}}
\clh$} \put(36,3){0}

 \put(6.6,2.1){$i_{\bm{z}}$} \put(23.6, 2.1){$\pi_{\bm{z}}$}

 \put(3,6.5){$N_{\bm{z}}$} \put(15,6.4){$P_{\bm{z}}$} \put(28,6.4){$Q_{\bm{z}}$}

 \put(-5, 9.5){0} \put(2,9.5){$\cls$}\put(12,9.5){$\clr \otimes \cle$} \put(27,9.5){$\clh$} \put(36,9.5){0}
 \put(5.6,10.5){$i$} \put(22.6, 10.5){$\pi$}

 \put(-4,10){ \vector(1,0){5}} \put(5.5,3.5){ \vector(1,0){3}} \put(4.5,10){ \vector(1,0){5}} \put(21.5,10)
 { \vector(1,0){4}} \put(22,3.5){ \vector(1,0){3}} \put(31,3.5){ \vector(1,0){4}}  \put(31,10){ \vector(1,0){4}}
 \put(2.4,9.2){ \vector(0,-1){5}} \put(14,9.2){ \vector(0,-1){5}} \put(27,9.2){ \vector(0,-1){5}}

 \end{picture}
 \end{center}

\NI which is commutative with exact rows for all $\bm{w}$ in
$\Omega$ (see \cite{DP}). Here $N_{\bm{z}}, P_{\bm{z}}$ and
$Q_{\bm{z}}$ are the quotient module maps. Since one can identify
$\clh/ I_{\bm{z}} \clh$ with $l^2_m$ and $(\clr \otimes
\cle)/I_{\bm{z}} (\clr \otimes \cle)$ with $\cle$, the kernel
functions of $\clh$ and $\clr \otimes \cle$ are given by $Q_{\bm{z}}
Q_{\bm{w}}^*$ and $P_{\bm{z}} P_{\bm{w}}^*$, respectively. Moreover,
since $Q_{\bm{w}} \pi = \pi_{\bm{w}} P_{\bm{w}}$ for all $\bm{w} \in
\Omega$, it follows that \[Q_{\bm{z}} \pi \pi^* Q_{\bm{w}} =
\pi_{\bm{z}} P_{\bm{z}} P_{\bm{w}}^* \pi_{\bm{w}}^*. \quad \quad
(\z, \bm{w} \in \Omega)\] Using the fact that $\pi \pi^* = I_{\clh}$
and $P_{\bm{z}} P_{\bm{w}}^* = k(\bm{z}, \bm{w}) \otimes I_{\cle}$,
one can now conclude that
$$Q_{\bm{z}} Q_{\bm{w}}^* = k(\bm{z}, \bm{w}) \pi_{\bm{z}}
\pi_{\bm{w}}^*. \quad \quad (\bm{z}, \bm{w} \in \Omega)$$

\NI Conversely, let the kernel function of the quasi-free Hilbert
module $\clh$ has the factorization
$$K_{\clh}(\bm{z}, \bm{w}) = k(\bm{z}, \bm{w}) \pi_{\bm{z}}
\pi_{\bm{w}}^*, \quad \quad (\bm{z}, \bm{w} \in \Omega)$$ for some
function $\pi : \Omega \raro \cll(\cle, l^2_m)$. Note that if the
function $\pi$ satisfies the above equality then it is holomorphic
on $\Omega$. Define a linear map $X : \clh \raro \clr \otimes \cle$
so that $$X Q_{\bm{z}}^* \eta = P_{\bm{z}}^* \pi_{\bm{z}}^* \eta.
\quad \quad (\eta \in l^2_m)$$  It then follows that

\[
\langle X(Q^*_{\bm{w}} \eta), X (Q^*_{\bm{z}} \zeta) \rangle =
\langle P_{\bm{w}}^* \pi_{\bm{w}}^* \eta, P^*_{\bm{z}}
\pi^*_{\bm{z}} \zeta \rangle = \langle \pi_{\bm{z}} P_{\bm{z}}
P^*_{\bm{w}} \pi^*_{\bm{w}} \eta, \zeta \rangle  = \langle
Q_{\bm{z}} Q^*_{\bm{w}} \eta, \zeta \rangle  = \langle Q^*_{\bm{w}}
\eta, Q^*_{\bm{z}} \zeta \rangle,
\]for all $\eta, \zeta \in l^2_m$. Therefore, since $\{Q_{\bm{z}}^*
\eta : \bm{z} \in \Omega, \eta \in l^2_m\}$ is a total set of
$\clh$, that $X$ extends to a bounded isometric operator. Moreover,
by the reproducing property of the kernel function, it follows that
$$M_{z_i}^* X (Q_{\bm{z}}^* \eta) = M_{z_i}^* P^*_{\bm{z}}
(\pi^*_{\bm{z}} \eta) = \overline{z}_i P^*_{\bm{z}} \pi^*_{\bm{z}}
\eta = \overline{z}_i X(Q^*_{\bm{z}} \eta) = XQ^*_{\bm{z}}
(\overline{z}_i \eta) = X M_{z_i}^* (Q_{\bm{z}}^* \eta),$$ for all
$1 \leq i \leq n$ and $\eta \in l^2_m$. Hence, $X \in \clb(\clh,
\clr \otimes \cle)$ is a co-module map. \qed

The following result is an application of the previous theorem.

\begin{Theorem}\label{TH2}
Let $\clh$ be a quasi-free Hilbert module of finite multiplicity and
$\clr$ be a reproducing kernel Hilbert module over $A(\Omega)$ (or
over $\mathbb{C}[\bm{z}]$). Let $k$ be the kernel function of
$\clr$. Then $\clr \otimes \cle$ is a dilation of $\clh$ for some
Hilbert space $\cle$ if and only if
$$K_{\clh}(\bm{z}, \bm{w}) = k(\bm{z}, \bm{w}) \tilde{K}(\bm{z},
\bm{w}), \quad \quad (\z, \w \in \Omega)$$ for some positive
definite kernel $\tilde{K}$ over $\Omega$. Moreover, if $k^{-1}$ is
defined, then the above conclusion is true if and only if $k^{-1}
K_{\clh}$ is a positive definite kernel.
\end{Theorem}
\NI \textsf{Proof.} The necessary part follows from the previous
theorem by setting $\tilde{K}(\bm{z}, \bm{w}) = \pi_{\bm{z}}
\pi_{\bm{w}}^*$. To prove the sufficiency part, let $K_{\clh} = k
\cdot \tilde{K}$ for some positive definite kernel $\tilde{K}$. We
let $\clh(\tilde{K})$ be the corresponding reproducing kernel
Hilbert space and set $\cle = \clh(\tilde{K})$. Let \[\pi_{\bm{z}} =
ev_{\bm{z}} \in \clb(\cle, l^2_m) \quad \quad (\z \in \Omega)\]be
the evaluation operator for the reproducing kernel Hilbert space
$\clh(\tilde{K})$. Then \[\tilde{K}(\z, \w) = \pi_{\z} \pi^*_{\w}.
\quad \quad (\z, \w \in \Omega)\] Consequently, by the previous
theorem it follows that $\clr \otimes \cle$ is a dilation of $\clh$.
\qed

Note that the reproducing kernel Hilbert space corresponding to the
kernel function $\tilde{K}$ is not necessarily a bounded module over
$A(\Omega)$ or even over $\mathbb{C}[\bm{z}]$. If it is a bounded
module, then one can identify $\clm$ canonically with the Hilbert
module tensor product, $\clr \otimes_{\mathbb{C}[\z]}
\clh(\tilde{K})$, which yields an explicit representation of the
co-isometry from the co-extension space $\clr \otimes
\clh(\tilde{K})$ to $\clm$.

\subsection{Hereditary functional calculus}
Let $p$ be a polynomial in the $2 n$ variables $\z = (z_1, \ldots,
z_n), \bar{\w} = (\bar{w}_1, \ldots, \bar{w}_n)$, where the $\z$
-variables all commute and the $\bar{\w}$-variables all commute with
no assumptions made about the relation of the $\z$ and $\bar{\w}$
variables. For any commuting $n$-tuple of operators $\mathbf T =
(T_1, \ldots , T_n)$, define the value of $p$ at $\mathbf T$ using
the hereditary functional calculus (following Agler \cite{AIEOT}):
$$
p(T, T^*) = \sum_{\bm{k}, \bm{l}} a_{\bm{k}, \bm{l}} T^{\bm{k}}
{T^*}^{\bm{l}},
$$
where $p(\z,\bar{\w}) = \mathop{\sum}_{\bm{k}, \bm{l}} a_{\bm{k},
\bm{l}} \z^{\bm{k}} \bar{\w}^{\bm{l}}$ and $\bm{k}, \bm{l} \in
\mathbb{N}^n$. Here, in the ``non-commutative polynomial'' $p(\z,
\bar{\w})$, the ``$\z$'s'' are all placed on the left, while the
``$\bar{\w}$'s'' are placed on the right.

Let $\clr \subseteq \clo(\Omega, \cle)$ be an $\cle$-valued
reproducing kernel Hilbert module over $\Omega$ for some Hilbert
space $\cle$ and $k$ be a positive definite kernel over $\Omega$.
Moreover, let \[k^{-1}(\bm{z}, \bm{w}) = \mathop{\sum}_{\bm{k},
\bm{l}} a_{\bm{k}, \bm{l}} \z^{\bm{k}} \bar{\w}^{\bm{l}},\] be a
polynomial in $\bm{z}$ and $\bar{\bm{w}}$. Therefore, for the module
multiplication operators on $\clr$ one gets
\[ k^{-1}({M}, {M}^*) = \mathop{\sum}_{\bm{k}, \bm{l}}
a_{\bm{k}, \bm{l}} {M}^{\bm{k}} {M}^{*\bm{l}}.\]

\begin{Proposition}\label{K-pos}
Let $\clr\subseteq \clo(\Omega, \cle)$ be a reproducing kernel
Hilbert module with kernel function $K_{\clr}$. Moreover, let $k$ be
a positive definite function defined on $\Omega$ and $k^{-1}(\bm{z},
\bm{w}) = \mathop{\sum}_{\bm{k}, \bm{l}} a_{\bm{k}, \bm{l}}
\z^{\bm{k}} \bar{\w}^{\bm{l}}$ be a polynomial in $\bm{z}$ and
$\bar{\bm{w}}$. Then  \[k^{-1}(M, M^*) \geq 0,\]if and only if
\[(\z, \w) \mapsto k^{-1}(\z,\w) K_{\clr}(\z,\w),\]is a positive
definite kernel on $\Omega$.
\end{Proposition}
\NI\textsf{Proof.} For each $\z, \w \in \Omega$ and $\eta, \zeta \in
\cle$, as a result of the preceding identity,
\[\begin{split} \langle k^{-1}({M}, {M}^*) K_{\clr}(\cdot, \w) \eta
, K_{\clr}(\cdot, \z) \zeta \rangle_{\clr} &= \langle
\Big(\mathop{\sum}_{\bm{k}, \bm{l}} a_{\bm{k}, \bm{l}} {M}^{\bm{k}}
{M}^{*\bm{l}}\Big) K_{\clr}(\cdot, \w) \eta , K_{\clr}(\cdot, \z) \zeta \rangle_{\clr} \\
&= \mathop{\sum}_{\bm{k}, \bm{l}} a_{\bm{k}, \bm{l}} \langle
{M}^{*\bm{l}} K_{\clr}(\cdot, \w) \eta , {M}^{*\bm{k}} K_{\clr}(\cdot, \z) \zeta \rangle_{\clr} \\
&= \mathop{\sum}_{\bm{k}, \bm{l}} a_{\bm{k}, \bm{l}} \bm{z}^{\bm{k}}
\bar{\bm{w}}^{\bm{l}} \langle
K_{\clr}(\cdot, \w) \eta , K_{\clr}(\cdot, \z) \zeta \rangle_{\clr} \\
&= k^{-1}(\z,\w) \langle K_{\clr}(\z,\w) \eta, \zeta
\rangle_{\cle}\\&= \langle k^{-1}(\z,\w) K_{\clr}(\z,\w) \eta, \zeta
\rangle_{\cle}.
\end{split}\]
Hence, for $\{\bm{z}_i\}_{i=1}^l \subseteq \Omega$ and
$\{\eta_i\}_{i=1}^l \subseteq l^2_m$ and $l \in \mathbb{N}$ it
follows that
\[\begin{split}
\langle k^{-1}({M}, {M}^*) (\sum_{i=1}^l K_{\clr}(\cdot, \bm{z}_i)
\eta_i), & \sum_{j=1}^l K_{\clr}(\cdot, \bm{z}_j) \eta_j
\rangle_{\clr} \\&  = \sum_{i, j = 1}^l \langle k^{-1}({M}, {M}^*)
(K_{\clr}(\cdot, \bm{z}_i) \eta_i), K_{\clr}(\cdot, \bm{z}_j) \eta_j
\rangle_{\clr} \\&  = \sum_{i, j = 1}^l \langle k^{-1}(\bm{z}_j,
\bm{z}_i) K_{\clr}(\bm{z}_j, \bm{z}_i) \eta_i, \eta_j \rangle_{\cle}
\\& = \sum_{i, j = 1}^l \langle (k^{-1} \circ K_{\clr})(\bm{z}_j,
\bm{z}_i) \eta_j, \eta_i\rangle_{\cle}.
\end{split}\]
Consequently, $k^{-1}(M, M^*) \geq 0$ if and only if $k^{-1}(\z,\w)
K_{\clm}(\z,\w)$ is a non-negative definite kernel. This completes
the proof. \qed

The following corollary is immediate.

\begin{Corollary}\label{TH3} Let $\clr\subseteq \clo(\Omega, \cle)$
be a reproducing kernel Hilbert module with kernel function
$K_{\clr}$. Moreover, let $k$ be a positive definite function
defined on $\Omega$ and $k^{-1}(\bm{z}, \bm{w}) =
\mathop{\sum}_{\bm{k}, \bm{l}} a_{\bm{k}, \bm{l}} \z^{\bm{k}}
\bar{\w}^{\bm{l}}$ be a polynomial in $\bm{z}$ and $\bar{\bm{w}}$.
Then  $k^{-1}(M, M^*) \geq 0$ if and only if $K_{\clr}$ factorizes
as
$$
K_{\clr}(\z,\w) = k(\z,\w) \tilde{K}(\z,\w),\quad \quad (\z,\w \in
\Omega)
$$
for some positive definite kernel $\tilde{K}$ on $\Omega$.
\end{Corollary}

The following dilation result is an application of Theorem \ref{TH2}
and Corollary \ref{TH3}.

\begin{Theorem} \label{dil}
Let $\clm$ be a quasi-free Hilbert module over $A(\mathbb{D}^n)$ of
multiplicity $m$ and $\clh_k$ be a reproducing kernel Hilbert module
over $A(\mathbb{D}^n)$. Moreover, let $k^{-1}(\bm{z}, \bm{w}) =
\mathop{\sum}_{\bm{k}, \bm{l}} a_{\bm{k}, \bm{l}} \z^{\bm{k}}
\bar{\w}^{\bm{l}}$ be a polynomial in $\bm{z}$ and $\bar{\bm{w}}$.
Then $\clh_k \otimes \clf$ is a dilation of $\clm$ for some Hilbert
space $\clf$ if and only if $k^{-1} ({M}, {M}^*) \geq 0$.
\end{Theorem}

It is the aim of the present consideration to investigate the issue
of uniqueness of the minimal isometric dilations of contractive
reproducing kernel Hilbert modules. The proof is based on operator
theory exploiting the fact that the co-ordinate multipliers define
doubly commuting isometries.

\begin{Theorem}\label{unique_H}
Let $\clh_k$ be a contractive reproducing kernel Hilbert module over
$A(\mathbb{D}^n)$. Then $\clh_k$ dilates to $H^2(\mathbb{D}^n)
\otimes \cle$ if and only if \, $\mathbb{S}^{-1}(M, M^*) \geq 0$ or,
equivalently, $\mathbb{S}^{-1} k \geq 0$. Moreover, if such dilation
exists, then the minimal one is unique.
\end{Theorem}
\NI\textsf{Proof.} By virtue of Theorem \ref{dil}, one only needs to
prove the uniqueness of the minimal dilation. Let $\Pi_i : \clh_k
\raro H^2(\mathbb{D}^n) \otimes \cle_i$ be minimal isometric
dilations of $\clh_k$, that is,
\[H^2(\mathbb{D}^n) \otimes \cle_i =
\overline{\mbox{span}}\{M_z^{\bm{k}} (\Pi_i \clh_k) : \bm{k} \in
\mathbb{N}^n\},\]for $i = 1, 2$. Define \[V : H^2(\mathbb{D}^n)
\otimes \cle_1 \raro H^2(\mathbb{D}^n) \otimes \cle_2,\]by
\[V(\sum_{|\bm{k}| \leq N} M_z^{\bm{k}} \Pi_1 f_{\bm{k}}) =
\sum_{|\bm{\alpha}| \leq N} M_z^{\bm{k}} \Pi_2 f_{\bm{k}},\]where
$f_{\bm{k}} \in \clh$ and $N \in \mathbb{N}$. Let $\bm{k}, \bm{l}
\in \mathbb{N}^n$ and define multi-indices $\tilde{\bm{k}}$ and
$\tilde{\bm{l}}$ so that
\[ {\tilde{k}}_i  = \left\{ \begin{array}{cc}
k_i - l_i & \mbox{for}\, k_i - l_i \geq 0  \\
0 & \mbox{for}\, k_i - l_i < 0 \end{array} \right.  \quad \mbox{and}
\quad \tilde{l}_i  = \left\{ \begin{array}{cc}
l_i - k_i & \mbox{for}\, l_i - k_i \geq 0  \\
0 & \mbox{for}\, l_i - k_i < 0 \end{array} \right.\]Note that $k_i -
l_i = \tilde{k}_i - \tilde{l}_i, \tilde{k}_i, \tilde{l}_i \geq 0$
and hence \[M_z^{*\bm{l}} M_z^{\bm{k}} = M_z^{*\tilde{\bm{l}}}
M_z^{\tilde{\bm{k}}} = M_z^{\tilde{\bm{k}}} M_z^{*\tilde{\bm{l}}}.\]
Therefore, for $i= 1, 2$, it follows that \[\langle M_z^{\bm{k}}
\Pi_i f_{\bm{k}}, M_z^{\bm{l}} \Pi_i f_{\bm{l}} \rangle = \langle
M_z^{*\bm{l}} M_z^{\bm{k}} \Pi_i f_{\bm{k}},  \Pi_i f_{\bm{l}}
\rangle = \langle M_z^{*\tilde{\bm{l}}} \Pi_i f_{\bm{k}},
M_z^{*\tilde{\bm{k}}} \Pi_i f_{\bm{l}} \rangle,\]and, since $\Pi_i$
is an co-module isometry, one gets \[\langle M_z^{\bm{k}} \Pi_i
f_{\bm{k}}, M_z^{\bm{l}} \Pi_i f_{\bm{l}} \rangle = \langle \Pi_i
M_z^{*\tilde{\bm{l}}} f_{\bm{k}}, \Pi_i M_z^{*\tilde{\bm{k}}}
f_{\bm{l}} \rangle = \langle M_z^{*\tilde{\bm{l}}} f_{\bm{k}},
M_z^{*\tilde{\bm{k}}} f_{\bm{l}} \rangle .\]Hence $V$ is
well-defined and isometric and \[V \Pi_1 = \Pi_2.\] Moreover, since
$$\{\sum_{|\bm{k}| \leq N} M_z^{\bm{k}} \Pi_i f_{\bm{k}}: f_{\bm{k}}
\in \clh, N \in \mathbb{N}\}$$ is a total subset of
$H^2(\mathbb{D}^n) \otimes \cle_i$ for $i=1, 2$, by minimality, $V$
is a unitary module map and hence $V = I_{H^2(\mathbb{D}^n)} \otimes
V_0$ for some unitary $V_0 \in \clb(\cle_1, \cle_2)$. Therefore, the
minimal dilations $\Pi_1$ and $\Pi_2$ are unitarily equivalent,
which concludes the proof. \qed

\begin{Corollary}\label{cor1}
If $\clh_k$ be a contractive reproducing kernel Hilbert space over
$A(\mathbb{D}^n)$. Then the Hardy module $H^2(\mathbb{D}^n) \otimes
\cle$ is a dilation of $\clh_k$ if and only if $\mathbb S_n^{-1}(M,
M^*) \geq 0$ or, equivalently, if and only if $\mathbb S_n^{-1} k
\geq 0$. Moreover, if an $H^2 (\mathbb{D}^n) \otimes {\mathcal E}$
dilation exists, then the minimal one is unique.
\end{Corollary}
\NI\textsf{Proof.} The necessary and sufficient part follows from
Theorem \ref{dil}. The uniqueness part follows from Theorem
\ref{unique_H}. \qed

The above proof will only work if the algebra is generated by
functions for which module multiplication  defines doubly commuting
isometric operators which happens for the Hardy module on the
polydisk. For a more general quasi-free Hilbert module $\clr$, the
maps $X^*_i$ identify anti-holomorphic sub-bundles of the bundle
$E_{\clr} \otimes \cle_i$, where $E_{\clr}$ is the Hermitian
holomorphic line bundle defined by $\clr$. To establish uniqueness,
some how one must extend this identification to the full bundles.
Equivalently, one has to identify the holomorphic quotient bundles
of $E_{\clr} \otimes \cle_1$, and $E_{\clr} \otimes \cle_2$ and must
some how lift it to the full bundles. At this point it is not even
obvious that the dimensions of $\cle_1$ and $\cle_2$ or the ranks of
the bundles are equal. This seems to be an interesting question.
Using results on exact sequences of bundles (cf. \cite{GH} and
\cite{W}), one can establish uniqueness if $\mbox{dim}\, \cle =
\mbox{rank}\, E_{\clh} + 1$.

\NI \textbf{Further results and comments:}
\begin{enumerate}

\item Most of the material in this section is based on the article \cite{DMS}.

\item In \cite{AJOT}, \cite{AIEOT}, \cite{At}, \cite{A-90}, \cite{A-92},
\cite{AE} and \cite{AEM}, Agler, Athavale, Ambrozie, Arazy, Englis
and Muller pointed out that the dilation theory and operator
positivity implemented by kernel functions are closely related to
each other.

\item Theorem \ref{unique_H} was proved by Douglas and Foias in \cite{DF} for the case of
multiplicity one. More precisely, let $\cls_1$ and $\cls_2$ be two
submodules of $H^2(\mathbb{D}^n)$. Then $H^2(\mathbb{D}^n)/ \cls_1
\cong H^2(\mathbb{D}^n)/ \cls_2$ if and only if $\cls_1 = \cls_2$.
This is a rigidity result concerning submodules of the Hardy module
(see Section \ref{R}).

\item Notice that any $n$-tuple of doubly commuting contractions on a
functional Hilbert space over $A(\mathbb{D}^n)$ satisfies the
hypothesis of Theorem \ref{dil}. Consequently, one can recover the
result of Sz.-Nagy and Foias (cf. \cite{NF}) in this situation. In
particular, $\clm^n = \clm \otimes \cdots \otimes \clm$ always
possesses a dilation to the Hardy module $H^2(\mathbb D^n)\otimes
\mathcal E$, where $\mathcal E$ is some Hilbert space, if $\clm$ is
contractive Hilbert module. The contractivity condition implies that
\[K(\z,\w) = (1-z_\ell \bar{w}_\ell)^{-1} Q_\ell(\z,\w), \quad \quad
(\z,\w \in \mathbb D^n)\] for some positive definite kernel $Q_\ell$
and for each $\ell = 1,2,\ldots , n$. Thus
\[K^n(\z,\w) = \mathbb S_n(\z,\w) Q(\z,\w), \quad \quad (\z,\w \in
\mathbb D^n)\] where $Q=\prod_{\ell=1}^n Q_\ell$. Thus the Hilbert
module $\clm^n$ corresponding to the positive definite kernel  $K^n$
is contractive and admits the kernel $\mathbb S_{\mathbb D^n}$ as a
factor, as shown above. This shows that $\clm^n$ has an isometric
co-extension to $H^2_{\mathcal Q}(\mathbb{D}^n)$, where $\mathcal Q$
is the reproducing kernel Hilbert space for the kernel $Q$.

\end{enumerate}

\section{Hardy module over polydisc}\label{HMOP}

This section begins by formulating a list of basic problems in
commutative algebra. Let $\clm$ be a module over $\mathbb{C}[z]$ and
$\clm^{\otimes n} := \clm \otimes_{\mathbb{C}} \cdots
\otimes_{\mathbb{C}} \clm$, the $n$-fold vector space tensor product
of $\clm$. Then $\clm^{\otimes n}$ is a module over $\mathbb{C}[z]
\otimes_{\mathbb{C}} \cdots \otimes_{\mathbb{C}} \mathbb{C}[z] \cong
\mathbb{C}[\z]$. Here the module action on $\clm^{\otimes n}$ is
given by
\[(p_1 \otimes \cdots \otimes p_n) \cdot (f_1 \otimes \cdots \otimes
f_n) \mapsto p_1 \cdot f_1 \otimes \cdots \otimes p_n \cdot
f_n,\]for all $\{p_i\}_{i=1}^n \subseteq \mathbb{C}[z]$ and
$\{f_i\}_{i=1}^n \in \clm_i$. Let $\{\clq_i\}_{i=1}^n$ be quotient
modules of $\clm$. Then
\begin{equation}\label{T-Q} \clq_1 \otimes_{\mathbb{C}} \cdots
\otimes_{\mathbb{C}} \clq_n,\end{equation} is a quotient module of
$\clm^{\otimes n}$.

On the other hand, let $\clq$ be a quotient module and $\cls$ a
submodule of $\clm_n$. One is naturally led to formulate the
following problems:

\NI (a) When is $\clq$ of the form (\ref{T-Q})?

\NI (b) When is $\clm/ \cls$ of the form (\ref{T-Q})?

Let now $\clm$ be the Hardy space $H^2(\mathbb{D})$, the Hilbert
space completion of $\mathbb{C}[z]$, and consider the analogous
problem. The purpose of this section is to provide a complete answer
to these questions when $\clm = H^2(\mathbb{D})$. In particular, a
quotient module $\clq$ of the Hardy module $H^2(\mathbb{D}^n) \cong
H^2(\mathbb{D}) \otimes \cdots \otimes H^2(\mathbb{D})$ is of the
form
\[\clq = \clq_1 \otimes \cdots \otimes \clq_n,\]
for $n$ quotient modules $\{\clq_i\}_{i=1}^n$ of $H^2(\mathbb{D})$
if and only if $\clq$ is doubly commuting.

A quotient module $\clq \subseteq H^2(\mathbb{D}^n)$ is said to be
\textit{doubly commuting} if \[C_{z_i} C_{z_j}^* = C_{z_j}^*
C_{z_i}. \quad \quad (1 \leq i < j \leq n)\] A submodule $\cls$ is
called \textit{co-doubly commuting} if $\cls^\perp \cong
H^2(\mathbb{D}^n)/\cls$ is doubly commuting quotient module.

\subsection{Submodules and Jordan blocks}

A closed subspace $\clq \subseteq H^2(\mathbb{D})$ is said to be a
\textit{Jordan block} of $H^2(\mathbb{D})$ if $\clq$ is a quotient
module and $\clq \neq H^2(\mathbb{D})$ (see \cite{NF1}, \cite{NF}).
By Beurling's theorem (see Corollary 5.4 in \cite{JS-HB}), a closed
subspace $\clq (\neq H^2(\mathbb{D}))$ is a quotient module of
$H^2(\mathbb{D})$ if and only if the submodule $\clq^{\perp}$ is
given by $\clq^{\perp} = \Theta H^2(\mathbb{D})$ for some inner
function $\Theta \in H^{\infty}(\mathbb{D})$. In other words, the
quotient modules and hence the Jordan blocks of $H^2(\mathbb{D})$
are precisely given by
\[\clq_\Theta : = H^2(\mathbb{D})/\Theta H^2(\mathbb{D}),\]for inner
functions $\Theta \in H^{\infty}(\mathbb{D})$. Thus on the level of
orthogonal projections, one gets\[P_{\clq_{\Theta}} =
I_{H^2(\mathbb{D})} - M_{\Theta} M^*_{\Theta} \quad \; \mbox{and} \;
\quad P_{\Theta H^2(\mathbb{D})} = M_{\Theta} M_{\Theta}^*.\]

The following lemma is a variation on the theme of the isometric
dilation theory of contractions.

\begin{Lemma}\label{Q-L}
Let $\clq$ be a quotient module of $H^2(\mathbb{D})$ and $\cll =
\mbox{ran} (I_{\clq} - C_z C^*_z) = \mbox{ran} (P_{\clq}
P_{\mathbb{C}} P_{\clq})$. Then $\clq = \mathop{\vee}_{l=0}^{\infty}
P_{\clq} M_z^l \cll$.
\end{Lemma}
\NI\textsf{Proof.} The result is trivial if $\clq = \{0\}$. Let
$\clq \neq \{0\}$, that is, $\clq^\perp$ is a proper submodule of
$H^2(\mathbb{D})$, or equivalently, $1 \notin \clq^\perp$. Notice
that
\[\mathop{\vee}_{l=0}^{\infty} P_{\clq} M_z^l \cll \subseteq \clq.\]
Let now \[f = \sum_{l=0}^\infty a_l z^l \in \clq,\] be such that $f
\perp \vee_{l=0}^{\infty} P_{\clq} M_z^l \cll$. It then follows that
$f \perp P_{\clq} M_z^l P_{\clq} P_{\mathbb{C}} \clq$, or
equivalently, $P_{\mathbb{C}} M_z^{*l} f \in {\clq}^{\perp}$ for all
$l \geq 0$. Since $P_{\mathbb{C}} M_z^{*l} f = a_l \in \mathbb{C}$
and $1 \notin {\clq}^{\perp}$, it follows that $a_l = P_{\mathbb{C}}
M_z^{*l} f = 0$ for all $l \geq 0$. Consequently, $f = 0$. This
concludes the proof. \qed

\subsection{Reducing submodules}

The following result gives a characterization of $M_{z_1}$-reducing
subspace of $H^2(\mathbb{D}^n)$.

\begin{Proposition}\label{reducing}
Let $n >1$ and $\cls$ be a closed subspace of $H^2(\mathbb{D}^n)$.
Then $\cls$ is a $(M_{z_2}, \ldots, M_{z_n})$-reducing subspace of
$H^2(\mathbb{D}^n)$ if and only if $\cls = \cls_1 \otimes
H^2(\mathbb{D}^{n-1})$ for some closed subspace $\cls_1$ of
$H^2(\mathbb{D})$.
\end{Proposition}

\NI\textsf{Proof.} Let $\cls$ be a $(M_{z_2}, \ldots,
M_{z_{n}})$-reducing closed subspace of $H^2(\mathbb{D}^n)$, that
is, $M_{z_i} P_{\cls} = P_{\cls} M_{z_i}$ for all $2 \leq i \leq n$.
Since
\[\begin{split} \sum_{\substack{0 \leq i_1 < \ldots < i_l \leq n, i_1, i_2 \neq 1}} (-1)^l
M_{z_{i_1}} \cdots M_{z_{i_l}} M^*_{{z}_{i_1}} \cdots
M^*_{{z}_{i_l}} & = (I_{H^2(\mathbb{D}^n)} - M_{z_2} M_{z_2}^*)
\cdots (I_{H^2(\mathbb{D}^2)} - M_{z_n} M_{z_n}^*) \\ & =
P_{H^2(\mathbb{D})} \otimes P_{\mathbb{C}} \otimes \cdots \otimes
P_{\mathbb{C}},\end{split}\] we have that $(P_{H^2(\mathbb{D})}
\otimes P_{\mathbb{C}} \otimes \cdots \otimes P_{\mathbb{C}})
P_{\cls} = P_{\cls} (P_{H^2(\mathbb{D})} \otimes P_{\mathbb{C}}
\otimes \cdots \otimes P_{\mathbb{C}})$. Therefore, $P_{\cls}
(P_{H^2(\mathbb{D})} \otimes P_{\mathbb{C}} \otimes \cdots \otimes
P_{\mathbb{C}})$ is an orthogonal projection and
\[P_{\cls} (P_{H^2(\mathbb{D})} \otimes P_{\mathbb{C}} \otimes
\cdots \otimes P_{\mathbb{C}}) = (P_{H^2(\mathbb{D})} \otimes
P_{\mathbb{C}} \otimes \cdots \otimes P_{\mathbb{C}}) P_{\cls} =
P_{\tilde{\cls_1}},\]where $\tilde{\cls_1} := (H^2(\mathbb{D})
\otimes \mathbb{C} \otimes \cdots \otimes \mathbb{C}) \cap \cls$.
Let $\tilde{\cls}_1 = \cls_1 \otimes \mathbb{C} \otimes \cdots
\otimes \mathbb{C}$, for some closed subspace $\cls_1$ of
$H^2(\mathbb{D})$. It then follows that
\[\cls = \overline{\mbox{span}} \{M_{z_2}^{l_2} \cdots M_{z_{n}}^{l_{n}}
\tilde{\cls}_1 : l_2, \ldots, l_{n} \in \mathbb{N}\} = \cls_1
\otimes H^2(\mathbb{D}^{n-1}).\]The converse part is immediate. This
concludes the proof of the proposition. \qed

It is well known that a closed subspace $\clm$ of
$H^2(\mathbb{D}^n)$ is $M_{z_1}$-reducing if and only if
\[\clm = H^2(\mathbb{D}) \otimes \cle,\]for some closed subspace
$\cle \subseteq H^2(\mathbb{D}^{n-1})$. The following key
proposition is a generalization of this fact.

\begin{Proposition}\label{1-reducing}
Let $\clq_1$ be a quotient module of $H^2(\mathbb{D})$ and $\clm$ be
a closed subspace of $\clq = \clq_1 \otimes H^2(\mathbb{D}^{n-1})$.
Then $\clm$ is a $P_{\clq} M_{z_1}|_{\clq}$-reducing subspace of
$\clq$ if and only if \[\clm = \clq_1 \otimes \cle,\]for some closed
subspace $\cle$ of $H^2(\mathbb{D}^{n-1})$.
\end{Proposition}

\NI\textsf{Proof.} Let  $\clm$ be a $P_{\clq}
M_{z_1}|_{\clq}$-reducing subspace of $\clq$. Then
\begin{equation}\label{M-Q}(P_{\clq} M_{z_1}|_{\clq}) P_{\clm} =
P_{\clm} (P_{\clq} M_{z_1}|_{\clq}),\end{equation}or equivalently,
$(P_{\clq_1} M_{z}|_{\clq_1} \otimes I_{H^2(\mathbb{D}^{n-1})})
P_{\clm} = P_{\clm} (P_{\clq_1} M_{z}|_{\clq_1} \otimes
I_{H^2(\mathbb{D}^{n-1})})$. Now \[I_{\clq} - (P_{\clq}
M_{z_1}|_{\clq}) (P_{\clq} M_{z_1}|_{\clq})^* = (P_{\clq_1}
P_{\mathbb{C}}|_{\clq_1}) \otimes
I_{H^2(\mathbb{D}^{n-1})}.\]Further (\ref{M-Q}) yields $P_{\clm}
((P_{\clq_1} P_{\mathbb{C}}|_{\clq_1}) \otimes
I_{H^2(\mathbb{D}^{n-1})}) = ((P_{\clq_1} P_{\mathbb{C}}|_{\clq_1})
\otimes I_{H^2(\mathbb{D}^{n-1})}) P_{\clm}$. Let \[\cll : = \clm
\cap \mbox{ran~} ((P_{\clq_1} P_{\mathbb{C}}|_{\clq_1}) \otimes
I_{H^2(\mathbb{D}^{n-1})}) = \clm \cap (\cll_1 \otimes
H^2(\mathbb{D}^{n-1})),\]where
\[\cll_1 = \mbox{ran~} (P_{\clq_1} P_{\mathbb{C}}|_{\clq_1})
\subseteq \clq_1.\]Since $\cll \subseteq \cll_1 \otimes
H^2(\mathbb{D}) \otimes \cdots \otimes H^2(\mathbb{D})$ and
$\mbox{dim~} \cll_1 = 1$ (otherwise, by Lemma \ref{Q-L} that $\cll_1
= \{0\}$ is equivalent to $\clq_1 = \{0\}$), it follows that $\cll =
\cll_1 \otimes \cle$ for some closed subspace $\cle \subseteq
H^2(\mathbb{D}^{n-1})$. Hence $P_{\clm} ((P_{\clq_1}
P_{\mathbb{C}}|_{\clq_1}) \otimes I_{H^2(\mathbb{D}^{n-1})}) =
P_{\cll} = P_{\cll_1 \otimes \cle}$. Claim: \[\clm =
\mathop{\vee}_{l = 0}^{\infty} P_{\clq} M_{z_1}^l \cll.\]Since
$\clm$ is $P_{\clq} M_{z_1}|_{\clq}$-reducing subspace and $\clm
\supseteq \cll$, it follows that $\clm \supseteq \mathop{\vee}_{l =
0}^{\infty} P_{\clq} M_{z_1}^l \cll$. To prove the reverse
inclusion, we let $f \in \clm$ and  $f = \sum_{\bm{k} \in
\mathbb{N}^n} a_{\bm{k}} \bm{z}^{\bm{k}}$, where $a_{\bm{k}} \in
\mathbb{C}$ for all $\bm{k} \in \mathbb{N}^n$. Then $f = P_{\clm}
P_{\clq} f = P_{\clm} P_{\clq} \sum_{\bm{k} \in \mathbb{N}^n}
a_{\bm{k}} \bm{z}^{\bm{k}}$. Observe now that for all $\bm{k} \in
\mathbb{N}^n$,
\[\begin{split} P_{\clm} P_{\clq} \bm{z}^{\bm{k}} & = P_{\clm} ((P_{\clq_1} z_1^{k_1}) (z_2^{k_2}
\cdots z_n^{k_n}))\\ & = P_{\clm} ((P_{\clq_1} M_{z_1}^{k_1}
P_{\clq_1} 1) (z_2^{k_2} \cdots z_n^{k_n}))\\& = P_{\clq}
M_{z_1}^{k_1} (P_{\clm} (P_{\clq_1} 1 \otimes z_2^{k_2} \cdots
z_n^{k_n})),\end{split}\]by (\ref{M-Q}), where the second equality
follows from $\langle z_1^{k_1}, f \rangle = \langle 1,
(M_{z_1}^{k_1})^* f \rangle = \langle P_{\clq_1} M_{z_1}^{k_1}
P_{\clq_1} 1, f \rangle$, for all $f \in \clq_1$. By the fact that
$P_{\clq_1} 1 \in \cll_1$ one gets $P_{\clm} (P_{\clq_1} 1 \otimes
z_2^{k_2} \cdots z_n^{k_n}) \in \cll$ and infer $P_{\clm} P_{\clq}
\bm{z}^{\bm{k}} \in \mathop{\vee}_{l = 0}^{\infty} P_{\clq_1}
M_{z_1}^l \cll$ for all $\bm{k} \in \mathbb{N}^n$. Therefore $f \in
\vee_{l = 0}^{\infty} P_{\clq} M_{z_1}^l \cll$ and hence $\clm =
\vee_{l=0}^\infty P_{\clq} M_{z_1}^l \cll$, and the claim follows.
Finally, $\cll = \cll_1 \otimes \cle$ yields \[\clm =
\mathop{\vee}_{l = 0}^{\infty} P_{\clq} M_{z_1}^l \cll =
(\mathop{\vee}_{l=0}^{\infty} P_{\clq_1} M_{z_1}^l \cll_1) \otimes
\cle,\]and therefore by Lemma \ref{Q-L}, $\clm = \clq_1 \otimes
\cle$. The converse part is trivial. This finishes the proof.  \qed

\subsection{Tensor product of Jordan blocks}

Let $\clq_1, \ldots, \clq_n$ be $n$ quotient modules of
$H^2(\mathbb{D})$. Then the module multiplication operators on $\clq
= \clq_1 \otimes \cdots \otimes \clq_n$ are given by
\[\{I_{\clq_1} \otimes \cdots \otimes {P_{\clq_i}
M_{z}|_{\clq_i}} \otimes \cdots \otimes I_{\clq_n}\}_{i=1}^n,\]that
is, $\clq$ is a doubly commuting quotient module. The following
theorem provides a converse statement.

\begin{Theorem}\label{dc-Q}
Let $\clq$ be a quotient module of $H^2(\mathbb{D}^n)$. Then $\clq$
is doubly commuting if and only if there exists quotient modules
$\clq_1, \ldots, \clq_n$ of $H^2(\mathbb{D})$ such that \[\clq =
\clq_1 \otimes \cdots \otimes \clq_n.\]
\end{Theorem}

\NI\textsf{Proof.} Let $\clq$ be a doubly commuting quotient module
of $H^2(\mathbb{D}^n)$. Set \[\tilde{\clq}_1 =
\overline{\mbox{span}} \{z_2^{l_2} \cdots z_n^{l_n} \clq : l_2,
\ldots, l_n \in \mathbb{N}\},\] a joint $(M_{z_2}, \ldots,
M_{z_n})$-reducing subspace of $H^2(\mathbb{D}^n)$. By Proposition
\ref{reducing}, it follows that
\[\tilde{\clq}_1 = \clq_1 \otimes H^2(\mathbb{D}^{n-1}),\]for some closed subspace
$\clq_1$ of $H^2(\mathbb{D})$. Since $\tilde{\clq}_1$ is
$M_{z_1}^*$-invariant subspace, that $\clq_1$ is a $M_z^*$-invariant
subspace of $H^2(\mathbb{D})$, that is, $\clq_1$ is a quotient
module of $H^2(\mathbb{D})$.

\NI Note that $\clq \subseteq \tilde{\clq}_1$. Claim: $\clq$ is a
$M_{z_1}^*|_{\tilde{\clq}_1}$-reducing subspace of $\tilde{\clq}_1$,
that is, \[P_{\clq} (M_{z_1}^*|_{\tilde{\clq}_1}) =
(M_{z_1}^*|_{\tilde{\clq}_1}) P_{\clq}.\]

\NI In order to prove the claim, first observe that $C_{z_1}^*
C_{z_i}^l = C_{z_i}^l C_{z_1}^*$ for all $l \geq 0$ and $2 \leq i
\leq n$, and hence
\[C_{z_1}^* C_{z_2}^{l_2} \cdots C_{z_n}^{l_n} = C_{z_2}^{l_2} \cdots C_{z_n}^{l_n}
C_{z_1}^*,\]for all $l_2, \ldots, l_n \geq 0$, that is, $M_{z_1}^*
P_{\clq} M_{z_2}^{l_2} \cdots M_{z_n}^{l_n} P_{\clq} = P_{\clq}
M_{z_2}^{l_2} \cdots M_{z_n}^{l_n} M_{z_1}^* P_{\clq}$ or, \[
M_{z_1}^* P_{\clq} M_{z_2}^{l_2} \cdots M_{z_n}^{l_n} P_{\clq} =
P_{\clq} M_{z_1}^* M_{z_2}^{l_2} \cdots M_{z_n}^{l_n} P_{\clq}.\]
From this it follows that for all $f \in \clq$ and $l_2, \ldots, l_n
\geq 0$,\[(P_{\clq} M_{z_1}^*|_{\tilde{\clq}_1}) (z_2^{l_2} \cdots
z_n^{l_n} f) = P_{\clq} M_{z_1}^* (z_2^{l_2} \cdots z_n^{l_n} f) =
(M_{z_1}^* P_{\clq}) (z_2^{l_2} \cdots z_n^{l_n} f).\]

\NI Also by $P_{\clq} \tilde{\clq}_1 \subseteq \tilde{\clq}_1$ one
gets
\[P_{\clq}
P_{\tilde{\clq}_1} = P_{\tilde{\clq}_1} P_{\clq}
P_{\tilde{\clq}_1}.\] This yields
\[\begin{split}(P_{\clq} M_{z_1}^*|_{\tilde{\clq}_1}) (z_2^{l_2} \cdots
z_n^{l_n} f) & = (M_{z_1}^* P_{\clq}) (z_2^{l_2} \cdots z_n^{l_n} f)
\\ & = M_{z_1}^* P_{\clq} P_{\tilde{\clq}_1} (z_2^{l_2} \cdots
z_n^{l_n} f) \\ & = (M_{z_1}^*|_{\tilde{\clq}_1} P_{\clq})
(z_2^{l_2} \cdots z_n^{l_n} f),\end{split}\]for all $f \in \clq$ and
$l_2, \ldots, l_n \geq 0$, and therefore
\[P_{\clq} (M_{z_1}^*|_{\tilde{\clq}_1}) = (M_{z_1}^*|_{\tilde{\clq}_1})
P_{\clq}.\]Hence $\clq$ is a $M_{z_1}^*|_{\tilde{\clq}_1}$-reducing
subspace of $\tilde{\clq}_1 = \clq_1 \otimes H^2(\mathbb{D}) \otimes
\cdots \otimes H^2(\mathbb{D})$. Now by Proposition
\ref{1-reducing}, there exists a closed subspace $\cle_1$ of
$H^2(\mathbb{D}^{n-1})$ such that
\[\clq = \clq_1 \otimes \cle_1.\]
Moreover, since\[\mathop{\vee}_{l = 0}^{\infty} z_1^l \clq =
\mathop{\vee}_{l=0}^{\infty} z_1^l (\clq_1 \otimes \cle_1) =
H^2(\mathbb{D}) \otimes \cle_1,\]and $\mathop{\vee}_{l = 0}^{\infty}
z_1^l \clq $ is a doubly commuting quotient module of
$H^2(\mathbb{D}^n)$, it follows that $\cle_1 \subseteq
H^2(\mathbb{D}^{n-1})$, a doubly commuting quotient module of
$H^2(\mathbb{D}^{n-1})$.

\NI By the same argument as above, we conclude that $\cle_1 = \clq_2
\otimes \cle_2$, for some doubly commuting quotient module of
$H^2(\mathbb{D}^{n-2})$. Continuing this process, we have $\clq =
\clq_1 \otimes \cdots \otimes \clq_n$, where $\clq_1, \ldots,
\clq_n$ are quotient modules of $H^2(\mathbb{D})$. This completes
the proof. \qed

As a corollary, one can easily derive the following fact concerning
Jordan blocks of $H^2(\mathbb{D}^n)$.

\begin{Corollary}\label{cor-JB}
Let $\clq$ be a closed subspace of $H^2(\mathbb{D}^n)$. Then $\clq$
is doubly commuting quotient module if and only if there exists
$\{\Theta_i\}_{i=1}^n \subseteq H^{\infty}(\mathbb{D})$ such that
each $\Theta_i$ is either inner or the zero function for all $1 \leq
i \leq n$ and \[\clq = \clq_{\Theta_1} \otimes \cdots \otimes
\clq_{\Theta_n}.\]
\end{Corollary}

\subsection{Beurling's representation}
The aim of this subsection is to relate the Hilbert tensor product
structure of the doubly commuting quotient modules to the Beurling
like representations of the corresponding co-doubly commuting
submodules.

The following piece of notation will be used in the rest of the
subsection.

\NI\textit{\textsf{Let $\Theta_i \in H^{\infty}(\mathbb{D})$ be a
given function indexed by $i \in \{1, \ldots, n\}$. In what follows,
$\tilde{\Theta}_i \in H^{\infty}(\mathbb{D}^n)$ will denote the
extension function defined by
\[\tilde{\Theta}_i(\bm{z}) = \Theta_i(z_i),\]for all $\bm{z} \in
\mathbb{D}^n$.}}

The reader is referred to Lemma 2.5 in \cite{JS2} for a proof of the
following lemma.

\begin{Lemma}\label{P-F} Let $\{P_i\}_{i=1}^n$ be a collection of commuting orthogonal
projections on a Hilbert space $\clh$. Then \[\cll :=
\mathop{\sum}_{i=1}^n \mbox{ran} P_i,\] is closed and the orthogonal
projection of $\clh$ onto $\cll$ is given by
\[\begin{split}P_{\cll} & = P_1 (I - P_2) \cdots (I - P_n) +
P_2 (I - P_3) \cdots (I - P_n) + \cdots + P_{n-1} (I - P_n) + P_n\\
& = P_n (I - P_{n-1}) \cdots (I - P_1) + P_{n-1} (I - P_{n-2})
\cdots (I - P_1) + \cdots + P_2 (I - P_1) +
P_1.\end{split}\]Moreover,
\[P_{\cll}
= I - \mathop{\prod}_{i=1}^n (I - P_i).\]
\end{Lemma}

The following provides an explicit correspondence between the doubly
commuting quotient modules and the co-doubly commuting submodules of
$H^2(\mathbb{D}^n)$.

\begin{Theorem}\label{10.7}
Let $\clq$ be a quotient module of $H^2(\mathbb{D}^n)$ and $\clq
\neq H^2(\mathbb{D}^n)$. Then $\clq$ is doubly commuting if and only
if there exists inner functions ${\Theta}_{i_j} \in
H^{\infty}(\mathbb{D})$ for $1 \leq i_1 < \ldots < i_m \leq n$ for
some integer $m \in \{1, \ldots, n\}$ such that
\[\clq = H^2(\mathbb{D}^n)/ [\tilde{\Theta}_{i_1} H^2(\mathbb{D}^n) + \cdots +
\tilde{\Theta}_{i_m} H^2(\mathbb{D}^n)],\]where
$\tilde{\Theta}_{i_j}(\bm{z}) = \Theta_{i_j}(z_{i_j})$ for all
$\bm{z} \in \mathbb{D}^n$.
\end{Theorem}
\NI\textsf{Proof.} The proof follows from Corollary \ref{cor-JB} and
Lemma \ref{P-F}. \qed

The conclusion of this subsection concerns the orthogonal projection
formulae for the co-doubly commuting submodules and the doubly
commuting quotient modules of $H^2(\mathbb{D}^n)$. It can be treated
as a co-doubly commuting submodules analogue of Beurling's theorem
on submodules of $H^2(\mathbb{D})$.

\begin{Corollary}\label{S-proj}
Let $\clq$ be a doubly commuting submodule of $H^2(\mathbb{D}^n)$.
Then there exists an integer $m \in \{1, \ldots, n\}$ and inner
functions $\{\Theta_{i_j}\}_{j=1}^m \subseteq
H^{\infty}(\mathbb{D})$ such that
\[\clq^{\perp} = \mathop{\sum}_{1 \leq i_1 < \ldots < i_m \leq n}
\tilde{\Theta}_{i_j} H^2(\mathbb{D}^n),\]where
$\tilde{\Theta_i}(\bm{z})= \Theta_{i_j}(z_{i_j})$ for all $\bm{z}
\in \mathbb{D}^n$. Moreover, \[P_{\clq} = I_{H^2(\mathbb{D}^n)} -
\mathop{\Pi}_{j=1}^m (I_{H^2(\mathbb{D}^n)} -
M_{\tilde{\Theta}_{i_j}} M_{\tilde{\Theta}_{i_j}}^*),\]and
\[P_{\clq^\perp} = \mathop{\Pi}_{j=1}^m
(I_{H^2(\mathbb{D}^n)} - M_{\tilde{\Theta}_{i_j}}
M_{\tilde{\Theta}_{i_j}}^*).\]
\end{Corollary}

This is an immediate consequence of Theorem \ref{10.7}. The reason
to state Theorem \ref{10.7} explicitly is its usefulness.

\NI \textbf{Further results and comments:}

\begin{enumerate}
\item An efficient solution to the algebraic problems, posed
in the introduction of this section, would likely have practical
applications.

\item The study of the doubly commuting quotient modules of  $H^2(\mathbb{D}^2)$ was
initiated by Douglas and Yang in \cite{DY} and \cite{DY1} (also see
\cite{BCL}). Later in \cite{INS} Izuchi, Nakazi and Seto obtained
the tensor product classification of doubly commuting quotient
modules of $H^2(\mathbb{D}^2)$.

\item The results of this section can be found in the papers \cite{JS1} and \cite{JS2}.
For the base case $n= 2$, they were obtained by Izuchi, Nakazi and
Seto \cite{IN}, \cite{INS}.

\item The tensor product representations of doubly commuting quotient modules of $H^2(\mathbb{D}^2)$ has deep and far reaching
applications to the general study of submodules and quotient modules
of the Hardy module $H^2(\mathbb{D}^2)$. See the work by K. J.
Izuchi, K. H. Izuchi and Y. Izuchi \cite{III1}, \cite{III} and Yang
\cite{Y-JOT05}, \cite{Y-JFA05}.

\item The techniques embodied in this section can be used to give
stronger results concerning doubly commuting quotient modules of a
large class of reproducing kernel Hilbert modules over
$\mathbb{D}^n$ including the weighted Bergman modules (see
\cite{CDS}).

\item Other related work concerning submodules and quotient modules of $H^2(\mathbb{D}^n)$
appears in Berger, Coburn and Lebow \cite{BCL}, Yang \cite{R-JFA01},
\cite{Y-JFA05}, Guo and Yang \cite{GY} and the book by Chen and Guo
\cite{C-G}.

\item In connection with Beurling representations for submodules
of $H^2(\mathbb{D}^n)$, we refer the reader to Cotlar and Sadosky
\cite{CoSa}.

\end{enumerate}

\section{Similarity to free Hilbert modules}\label{SFHM}

This section begins by describing the notion of "split short exact
sequence"from commutative algebra. Let $M_1$ and $M_2$ be modules
over a ring $R$. Then $M_1 \oplus M_2$, module direct sum of $M_1$
and $M_2$, yields the short exact sequence \[0 \longrightarrow M_1
\stackrel{i}{\longrightarrow} M_1 \oplus M_2 \stackrel{\pi}
\longrightarrow M_2 \longrightarrow 0,\]where $i$ is the embedding,
$m_1 \mapsto (m_1, 0)$ and $\pi$ is the projection, $(m_1, m_2)
\mapsto m_2$ for all $m_1 \in M_1$ and $m_2 \in M_2$. A short exact
sequence of modules  \[0 \longrightarrow M_1
\stackrel{\varphi_1}{\longrightarrow} M \stackrel{\varphi_2}
\longrightarrow M_2 \longrightarrow 0,,\]is called \textit{split
exact sequence} if there exists a module isomorphism $\varphi : M
\raro M_1 \oplus M_2$ such that the diagram

\[\begin{CD}
0 @>>>M_1 @> \varphi_1 >> M @> \varphi_2>>M_2 @ >>> 0\\
@. @V I_{M_1} VV@V \varphi VV @V I_{M_2}VV\\
0 @>>>M_1 @> i >> M_1\oplus M_2 @> \pi>>M_2 @ >>> 0
\end{CD}\]commutes. It is well known that a short exact sequence of modules
\[0 \longrightarrow M_1
\stackrel{\varphi_1}{\longrightarrow} M \stackrel{\varphi_2}
\longrightarrow M_2 \longrightarrow 0,\] splits if and only if
$\varphi_2$ has a right inverse, if and only if $\varphi_1$ has a
left inverse.

Now let $\clr \subseteq \clo(\Omega, \mathbb{C})$ be a reproducing
kernel Hilbert module. In the rest of this section focus will be on
the quotient module $\clq_\Theta$ of $\clr \otimes \cle_*$ given by
the exact sequence of Hilbert modules\[\cdots \longrightarrow \clr
\otimes \cle \stackrel{M_\Theta}{\longrightarrow} \clr \otimes
\cle_* \stackrel{\pi_\Theta} \longrightarrow \clq_\Theta
\longrightarrow 0,\]where $\Theta \in \clm_\clr(\cle, \cle_*)$ is a
multiplier. In other words, \[\clq_\Theta = (\clr \otimes
\cle_*)/\mbox{ran} M_\Theta.\] The exact sequence is called
\textit{split} if $\pi_\Theta$ has a module right inverse, that is,
if there exists a module map $\sigma_\Theta : \clq_\Theta \raro \clr
\otimes \cle_*$ such that
\[\pi_\Theta \sigma_\Theta = I_{\clr \otimes \cle_*}.\]

\subsection{Complemented submodules}

This subsection provides a direct result concerning splitting of
Hilbert modules, which involves a mixture of operator theory and
algebra.
\begin{Theorem}\label{1}
Let $\clr \subseteq \clo(\Omega, \mathbb{C})$ be a reproducing
kernel Hilbert module and $\Theta \in \clm_{\clb(\cle,
\cle_*)}(\clr)$ be a multiplier for Hilbert spaces $\cle$ and
$\cle_*$ such that $\mbox{ran}\,M_{\Theta}$ is closed. Then
$\mbox{ran}\,M_{\Theta}$ is complemented in $\clr \otimes \cle_*$ if
and only if $\pi_{\Theta}$ is right invertible, that is, there
exists a module map $\sigma_{\Theta} : \clq_{\Theta} \raro \clr
\otimes \cle_*$ such that $\pi_{\Theta} \sigma_\Theta = I_{\clr
\otimes \cle_*}$.
\end{Theorem}
\NI \textsf{Proof.} Let $\cls$ be a submodule of $\clr \otimes
\cle_*$ such that $\clr \otimes \cle_* = \mbox{ran}\,M_{\Theta}
\stackrel{\cdot} + \cls$. Then \[Y = \pi_{\Theta}|_{\cls} : \cls
\raro \clq_\Theta\] is one-to-one. Also for $f \in \clq_\Theta$ and
$f = f_1 + f_2$ with $f_1 \in \mbox{ran} M_\Theta$ and $f_2 \in
\cls$ one gets
\[\pi_\Theta f_2 = \pi_\Theta(f - f_1) = f - \pi_\Theta f_1 = f.\]Consequently,
$Y = \pi_\Theta|_\cls$ is onto. Hence \[Y^{-1} : (\clr \otimes
\cle_*)\,/\,\mbox{ran}\,M_{\Theta} \raro \cls,\] is bounded by the
closed graph theorem. Let \[\sigma_{\Theta} = i ~Y^{-1},\] where $i
: \cls \raro \clr \otimes \cle_*$ is the inclusion map. Again for $f
= f_1 + f_2 \in \clq_\Theta$ with $f_1 \in \mbox{ran} M_\Theta$ and
$f_2 \in \cls$ one gets \[\pi_\Theta \sigma_\Theta f = \pi_\Theta (i
Y^{-1}) (f_1 + f_2) = \pi_\Theta i f_2 = \pi_\Theta f_2 =
\pi_\Theta(f - f_1) = \pi_\Theta f = f.\]Therefore, $\sigma_\Theta$
is a right inverse for $\pi_{\Theta}$. Also that $\sigma_\Theta$ is
a module map follows from the fact that $Y$ is a module map.

\NI Conversely, let $\sigma_{\Theta} : \clq_{\Theta} \raro \clr
\otimes \cle_*$ be a module map which is a right inverse of
$\pi_{\Theta}$. Then $\sigma_{\Theta} \pi_{\Theta}$ is an idempotent
on $\clr \otimes \cle_*$ such that $\cls =
\mbox{ran}\,\sigma_{\Theta} \pi_{\Theta}$ is a complementary
submodule for the closed submodule $\mbox{ran}\,M_{\Theta}$ in $\clr
\otimes \cle_*$. \qed

Examples in the case $n = 1$ show that the existence of a right
inverse  for $\pi_{\Theta}$ does not imply that $\mbox{ran}\,
M_{\Theta}$ is closed. However, if $\Theta \in \clm_{\clb(\cle,
\cle_*)}(\clr)$ and $\mbox{ran}\, M_{\Theta}$ is complemented in
$\clr \otimes \cle_*$, then $\mbox{ran}\, M_{\Theta}$ is closed.

\subsection{Lifting and range-inclusion theorems}

This subsection concerned with the study of lifting and
Drury-Arveson module (see \cite{AAM} or \cite{JS-HB}). The
\textit{commutant lifting theorem} will be used to extend some
algebraic results for the case of quotient modules of the
Drury-Arveson module. The commutant lifting theorem for the
Drury-Arveson module is due to Ball, Trent and Vinnikov, Theorem 5.1
in \cite{BTV}.

\begin{Theorem}\label{CLT}
Let $\cln$ and $\cln_*$ be quotient modules of $H^2_n \otimes \cle$
and $H^2_n \otimes \cle_*$ for some Hilbert spaces $\cle$ and
$\cle_*$, respectively. If $X : \cln \raro \cln_*$ is a bounded
module map, that is,   $$X P_{\cln} (M_{z_i} \otimes
I_{\cle})|_{\cln} = P_{\cln_*} (M_{z_i} \otimes
I_{\cle_*})|_{\cln_*} X,$$ for $i = 1, \ldots, n$, then there exists
a multiplier $\Phi \in \clm_{\clb(\cle, \cle_*)}(H^2_n)$ such that
$\|X\| = \|M_{\Phi}\|$ and  $P_{\cln_*} M_{\Phi} = X$.
\end{Theorem}

In the language of Hilbert modules, one has the following
commutative diagram
\[\begin{CD}H^2_n \otimes \cle @>M_{\Phi} >> H^2_n \otimes
\cle_*\\ @ V\pi_{\cln}VV @V\pi_{\cln_*}VV\\\cln @>X>> \cln_*
\end{CD}\]where $\pi_{\cln}$ and $\pi_{\cln_*}$ are the
quotient maps.

As one knows, by considering the $n= 1$ case, there is more than one
multiplier $\Theta \in \clm_{\clb(\cle, \cle_*)}(H^2_n)$ for Hilbert
spaces $\cle$ and $\cle_*$ with the same range and thus yielding the
same quotient. Things are even more complicated for $n>1$. However,
the following result using the commutant lifting theorem introduces
some order.
\begin{Theorem}\label{incl}
Let $\Theta \in \clm_{\clb(\cle, \cle_*)}(H^2_n)$ be a multiplier
with closed range for Hilbert spaces $\cle$ and $\cle_*$ and $\Phi
\in \clm_{\clb(\clf, \cle_*)}(H^2_n)$ for some Hilbert space $\clf$.
Then
\[\mbox{ran~} M_{\Phi} \subseteq \mbox{ran~}M_{\Theta},\]if and only if
\[\Phi = \Theta \Psi,\]for some multiplier $\Psi \in \clm_{\clb(\clf,
\cle)}(H^2_n)$.
\end{Theorem}
\NI \textsf{Proof.} If $\Psi \in \clm(\clf, \cle)$ such that $\Phi =
\Theta \Psi$, then $M_{\Phi} = M_{\Theta} M_{\Psi}$ and hence
\[\mbox{ran~} M_{\Phi} = \mbox{ran~} M_{\Theta} M_{\Psi} \subseteq
\mbox{ran~}M_{\Theta}.\] Suppose $\mbox{ran~} M_{\Phi} \subseteq
\mbox{ran~}M_{\Theta}$. Consider the module map $\hat{M_{\Theta}} :
(H^2_n \otimes \cle)~/~\mbox{ker~}M_{\Theta} \longrightarrow
\mbox{ran~}M_{\Theta}$ defined by
\[\hat{M_{\Theta}} \gamma_{\Theta} = M_{\Theta},\]where $\gamma_{\Theta}
: H^2_n \otimes \cle \longrightarrow (H^2_n \otimes
\cle)/~\mbox{ker~}M_{\Theta}$ is the quotient module map. Since
$\mbox{ran~} M_{\Theta}$ is closed that $\hat{M}_\Theta$ is
invertible.  Then \[\hat{X}: = \hat{M_{\Theta}}^{-1} :
\mbox{~ran~}M_{\Theta} \raro (H^2_n \otimes
\cle)/\mbox{~ker~}M_{\Theta}\] is bounded by the closed graph
theorem and so is \[\hat{X} M_{\Phi} : H^2_n \otimes \clf \raro
(H^2_n \otimes \cle)\,/\,\mbox{ker~} M_{\Theta}.\] By Theorem
\ref{CLT} there exists a multiplier $\Psi \in \clm_{\clb(\clf,
\cle)}(H^2_n)$ so that
$$\gamma_{\Theta} M_{\Psi} = \hat{X} M_{\Phi},$$ and hence
\[M_{\Theta} M_{\Psi} = (\hat{M_{\Theta}} \gamma_{\Theta}) M_{\Psi}
= \hat{M_{\Theta}} (\hat{X} M_{\Phi}) =  M_{\Phi},\] or $\Phi =
\Theta \Psi$ which completes the proof. \qed

\subsection{Regular inverse and similarity problem}
The purpose of this subsection is to establish an equivalent
condition which will allow one to tell when the range of a
multiplier will be complemented.

A multiplier $\Theta \in \clm_{\clb(\cle, \cle_*)}(H^2_n)$ is said
to have a regular inverse if there exists $\Psi \in
\clm_{\clb(\cle_*, \cle)}(H^2_n)$ such that $${\Theta}(\z)
{\Psi}(\z) {\Theta}(\z) = {\Theta}(\z). \quad \quad (\z \in
\mathbb{B}^n)$$

\begin{Theorem}\label{P1}
Let $\Theta \in \clm_{\clb(\cle, \cle_*)}(H^2_n)$. Then $\Theta$
admits a regular inverse if and only if $\mbox{ran}\, M_{\Theta}$ is
complemented in $H^2_n \otimes \cle_*$, or $$H^2_n \otimes \cle_* =
\mbox{ran}\, M_{\Theta} \stackrel{\cdot}+ \cls,$$ for some submodule
$\cls$ of $H^2_n \otimes \cle_*$.
\end{Theorem}
\NI \textsf{Proof.} If $H^2_n \otimes \cle_* = \mbox{ran}\,
M_{\Theta} \stackrel{\cdot}+ \cls$ for some (closed) submodule
$\cls$, then $\mbox{ran~}M_{\Theta}$ is closed. Consider the module
map $\hat{M_{\Theta}} : (H^2_n \otimes \cle)\,/ \,\mbox{ker}\,
M_{\Theta} \longrightarrow (H^2_n \otimes \cle_*)\,/\, \cls$ defined
by $$\hat{M_{\Theta}} \gamma_{\Theta} = \pi_{\cls} M_{\Theta},$$
where $\gamma_{\Theta} : H^2_n \otimes \cle \raro (H^2_n \otimes
\cle)\,/\,\mbox{ker}\, M_{\Theta}$ and $\pi_{\cls} : H^2_n \otimes
\cle_* \raro (H^2_n \otimes \cle_*)\,/\,\cls$ are quotient maps.
This map is one-to-one and onto and thus has a bounded inverse
$\hat{X} = \hat{M_{\Theta}}^{-1} : (H^2_n \otimes \cle_*)\,/ \,\cls
\raro (H^2_n \otimes \cle)\,/\, \mbox{ker~}M_{\Theta}$ by the closed
graph theorem. Since $\hat{X}$ satisfies the hypotheses of the
commutant lifting theorem, there exists $\Psi \in \clm_{\clb(\cle_*,
\cle)}(H^2_n)$ such that $\gamma_{\Theta} M_{\Psi} = \hat{X}
\pi_{\cls}$. Further, $\hat{M_{\Theta}} \gamma_{\Theta} = \pi_{\cls}
M_{\Theta}$ yields
\[\pi_{\cls} M_{\Theta} M_{\Psi} = \hat{M_{\Theta}} \gamma_{\Theta}
M_{\Psi} = \hat{M_{\Theta}} \hat{X} \pi_{\cls} = \pi_{\cls},\] and
therefore,\[\pi_{\cls}(M_{\Theta} M_{\Psi} M_{\Theta} - M_{\Theta})
= 0.\] Since $\pi_{\cls}$ is one-to-one on $\mbox{ran}\,
M_{\Theta}$, it follows that $M_{\Theta} M_{\Psi} M_{\Theta} =
M_{\Theta}$.

\NI Now suppose there exists $\Psi \in \clm_{\clb(\cle_*,
\cle)}(H^2_n)$ such that $M_{\Theta} M_{\Psi} M_{\Theta} =
M_{\Theta}$. This implies  that \[(M_{\Theta} M_{\Psi})^2 =
M_{\Theta} M_{\Psi},\] and hence $M_{\Theta} M_{\Psi}$ is an
idempotent. From the equality $M_{\Theta} M_{\Psi} M_{\Theta} =
M_{\Theta}$ one obtains both that $\mbox{ran}\,M_{\Theta} M_{\Psi}$
contains $\mbox{ran}\,M_{\Theta}$ and that $\mbox{ran}\, M_{\Theta}
M_{\Psi}$ is contained in $\mbox{ran}\, M_{\Theta}$. Therefore,
$$\mbox{ran}\, M_{\Theta} M_{\Psi} = \mbox{ran}\, M_{\Theta},$$ and
$$\cls = \mbox{ran}\,(I - M_{\Theta} M_{\Psi}),$$ is a complementary
submodule of $\mbox{ran~}M_{\Theta}$ in $H^2_n \otimes \cle_*$. \qed

\begin{Corollary}\label{TH1}
Assume $\Theta \in \clm_{\clb(\cle, \cle_*)}(H^2_n)$ for Hilbert
spaces $\cle$ and $\cle_*$ such that $\mbox{ran}\, M_{\Theta}$ is
closed and $\clh_{\Theta}$ is defined by \[H^2_n \otimes \cle
\stackrel{M_{\Theta}}\longrightarrow H^2_n \otimes \cle_*
\longrightarrow \clh_{\Theta} \longrightarrow 0.\] If
$\clh_{\Theta}$ is similar to $H^2_n \otimes \clf$ for some Hilbert
space $\clf$, then the sequence splits.
\end{Corollary}
\NI \textsf{Proof.} First, assume that there exists an invertible
module map $X : H^2_n \otimes \clf \raro \clh_{\Theta}$, and let
$\Phi \in \clm_{\clb(\clf, \cle_*)}(H^2_n)$ be defined by the
commutant lifting theorem, Theorem \ref{CLT}, so that \[\pi_{\Theta}
M_{\Phi} = X,\] where $\pi_{\Theta} : H^2_n \otimes \cle_* \raro
(H^2_n \otimes \cle_*)\,/\, \mbox{ran~}M_{\Theta}$ is the quotient
map. Since $X$ is invertible one gets $$H^2_n \otimes \cle_* =
\mbox{ran~} M_{\Phi} \stackrel{\bm{.}}{+}  \mbox{ran} M_{\Theta}.$$
Thus $\mbox{ran~}M_{\Theta}$ is complemented and hence it follows
from Theorem \ref{1} that the sequence splits. \qed

Finally, the following weaker converse to Corollary \ref{TH1} always
holds.

\begin{Corollary}
Let $\Theta \in \clm_{\clb(\cle, \cle_*)}(H^2_n)$ for Hilbert spaces
$\cle$ and $\cle_*$, and set $\clh_{\Theta} = (H^2_n \otimes
\cle_*)/ \mbox{~clos~}[\mbox{ran~}M_{\Theta}]$. Then the following
statements are equivalent:

\NI (i) $\Theta$ is left invertible, that is, there exists $\Psi \in
\clm_{\clb(\cle_*, \cle)}(H^2_n)$ such that $\Psi \Theta =
I_{\cle}$.

\NI (ii) $\mbox{ran~}M_{\Theta}$ is closed, $\mbox{ker~}M_{\Theta} =
\{0\}$ and $\clh_{\Theta}$ is similar to a complemented submodule
$\cls$ of $H^2_n \otimes \cle_*$.
\end{Corollary}

\NI\textsf{Proof.} If (i) holds, then $\mbox{ran~}M_{\Theta}$ is
closed and $\mbox{ker~}M_{\Theta} = \{0\}$. Further, $M_{\Theta}
M_{\Psi}$ is an idempotent on $H^2_n \otimes \cle_*$ such that
$\mbox{ran~} M_{\Theta} M_{\Psi} = \mbox{ran~}M_{\Theta}$ and
$\clh_{\Theta}$ is isomorphic to \[\cls = \mbox{ran~}(I - M_{\Theta}
M_{\Psi}) \subseteq H^2_n \otimes \cle_*,\] and \[H^2_n \otimes
\cle_* = \mbox{ran~}M_{\Theta} \stackrel{\bm{.}}{+} \cls,\] so
$\cls$ is complemented.

\NI Now assume that (ii) holds and there exists an isomorphism $X :
\clh_{\Theta} \raro \cls \subseteq H^2_n \otimes \cle_*$, where
$\cls$ is a complemented submodule of $H^2_n \otimes \cle_*$. Then
\[Y = X \pi_{\Theta} : H^2_n \otimes \cle_* \raro H^2_n \otimes
\cle_*,\] is a module map and hence there exists a multiplier $\Xi
\in \clm_{\clb(\cle_*, \cle_*)}(H^2_n)$ so that $Y = M_{\Xi}$. Since
$X$ is invertible, $\mbox{ran~}M_{\Xi} = \cls$, which is
complemented by assumption, and hence by Theorem \ref{P1} there
exists $\Psi \in \clm_{\clb(\cle_*, \cle_*)}(H^2_n)$ such that
$M_{\Xi} = M_{\Xi} M_{\Psi} M_{\Xi}$ or $M_{\Xi}(I - M_{\Psi}
M_{\Xi}) = 0$. Therefore,
$$\mbox{ran~}(I - M_{\Psi} M_{\Xi}) \subseteq \mbox{ker~}
M_{\Xi} = \mbox{ker~}Y = \mbox{ker~}\pi_{\Theta} =
\mbox{ran~}M_{\Theta}.$$ Applying Theorem \ref{incl}, we obtain
$\Phi \in \clm_{\clb(\cle_*, \cle)}(H^2_n)$ so that $$I - M_{\Psi}
M_{\Xi} = M_{\Theta} M_{\Phi}.$$ Thus using $M_{\Xi} M_{\Theta} = 0$
we see that $M_{\Theta} M_{\Phi} M_{\Theta} = (I - M_{\Psi} M_{\Xi})
M_{\Theta} = M_{\Theta}$, or, $M_{\Theta} = M_{\Theta} M_{\Phi}
M_{\Theta}$. Since $\mbox{ker~}M_{\Theta} = \{0\}$, we have
$M_{\Phi} M_{\Theta} = I_{H^2_n \otimes \cle}$, which completes the
proof. \qed

Theorem \ref{P1} and Corollary \ref{TH1} yields the main result of
the present subsection.

\begin{Corollary}\label{main-cor}
Given $\Theta \in \clm_{\clb(\cle, \cle_*)}(H^2_n)$ for Hilbert
spaces $\cle$ and $\cle_*$ such that $\mbox{ran~}M_{\Theta}$ is
closed, consider the quotient Hilbert module $\clh_{\Theta}$ given
by \[H^2_n \otimes \cle \stackrel{M_{\Theta}}\longrightarrow H^2_n
\otimes \cle_* \longrightarrow \clh_{\Theta} \longrightarrow 0.\] If
$\clh_{\Theta}$ is similar to $H^2_n \otimes \clf$ for some Hilbert
space $\clf$, then $\Theta$ has a regular inverse.
\end{Corollary}

\NI \textbf{Further results and comments:}

\begin{enumerate}

\item Most of the material in this section can be found in \cite{DFS}.

\item It is well known that (see \cite{NF}) a contractive Hilbert module
$\clh$ over $A(\mathbb{D})$ is similar to a unilateral shift if and
only if its characteristic function $\Theta_{\clh}$ has a left
inverse. Various approaches to this result have been given but the
present one uses the commutant lifting theorem and, implicitly, the
Beurling-Lax-Halmos theorem. In particular, the proof does not
involve, at least explicitly, the geometry of the dilation space for
the contraction.

\item For $n = 1$, Corollary \ref{main-cor} yields a more general result
concerning similarity of contractive Hilbert modules over
$A(\mathbb{D})$:
\begin{Theorem}
Let $\cle$ and $\cle_*$ be Hilbert spaces and $\Theta \in
H^\infty_{\clb(\cle, \cle_*)}(\mathbb{D})$ be a bounded analytic
function such that $\mbox{ker~} M_{\Theta} = \{0\}$ and
$\mbox{ran~}M_{\Theta}$ is closed. Then the quotient module
$\clh_{\Theta}$ given by \[0 \longrightarrow H^2_{\cle}(\mathbb{D})
\stackrel{M_{\Theta}}\longrightarrow H^2_{\cle_*}(\mathbb{D})
\longrightarrow \clh_{\Theta} \longrightarrow 0,\]is similar to
$H^2_{\clf}(\mathbb{D})$ for some Hilbert space $\clf$ if and only
if $\Theta \Psi \Theta = \Theta$ for some $\Psi \in
H^\infty_{\clb(\cle_*, \cle)}(\mathbb{D})$.
\end{Theorem}

\item If $\clh$ is a Hilbert module over
$\mathbb{C}[\z]$, Corollary \ref{main-cor} remains true under the
assumption that the analogue of the commutant lifting theorem holds
for the class of Hilbert modules under consideration. In particular,
Corollary \ref{main-cor} can be generalized to any other reproducing
kernel Hilbert modules where the kernel is given by a complete
Nevanlinna-Pick kernel (see \cite{AM-JFA}).

\item For other results concerning similarity in both commutative
and noncommutative setup see Popescu \cite{GP11a}.

\end{enumerate}

\section{Generalized canonical models and similarity}\label{GCMS}

This section describes conditions for certain quotient Hilbert
modules to be similar to the reproducing kernel Hilbert modules from
which they are constructed.

In particular, it is shown that the similarity criterion for a
certain class of quotient Hilbert modules is independent of the
choice of the basic Hilbert module ``building blocks" as in the
isomorphism case, so long as the multiplier algebras are the same.

\subsection{Corona pairs}
This subsection begins with the case in which the existence of a
left inverse for the multiplier depends only on a positive answer to
the corona problem for the domain.

\begin{Theorem}\label{sim1}
Let $\clr \subseteq \clo(\Omega, \mathbb{C})$ be a reproducing
kernel Hilbert module over $\mathbb{C}[\z]$. Assume that $\theta_1,
\theta_2, \psi_1, \psi_2 $ are in $\clm(\clr)$ and that $\theta_1
\psi_1 + \theta_2 \psi_2 =1$. Then the quotient Hilbert module
$\clr_{\Theta}= (\clr \otimes \mathbb{C}^2)/ M_{\Theta} \clr$ is
similar to $\clr$, where $M_{\Theta} f = \theta_1 f \otimes e_1 +
\theta_2 f \otimes e_2 \in \clh \otimes \mathbb{C}^2$ and $f \in
\clr$, with $\{e_1, e_2\}$ the standard orthonormal basis for
$\mathbb{C}^2$.
\end{Theorem}
\NI\textsf{Proof.} Let $R_{\Psi} : \clr \oplus \clr \raro \clr$ be
the bounded module map defined by $R_{\Psi} (f \oplus g) = \psi_1 f
+ \psi_2 g$ for $f, g \in \clr$. Note that
$$R_{\Psi}M_{\Theta}=I_{\clh},$$ or that $R_{\Psi}$ is a left
inverse for $M_\Theta$. Then for any $f \oplus g \in \clr \oplus
\clr$, one gets $$f \oplus g = (M_{\Theta} R_{\Psi}(f \oplus g)) + (
f \oplus g - M_{\Theta} R_{\Psi}(f \oplus g)),$$ with $M_{\Theta}
R_{\Psi}(f \oplus g) \in \mbox{ran}\, M_{\Theta}$ and $f \oplus g -
M_{\Theta} R_{\Psi}(f \oplus g) \in \mbox{ker}\, R_{\Psi}$. This
decomposition, along with $$\mbox{ran}\, M_{\Theta} \cap
\mbox{ker}\, R_{\Psi} = \{0\}$$ implies that $$\clr \oplus \clr =
\mbox{ran}\, M_{\Theta} \stackrel{\cd}+ \mbox{ker}\, R_{\Psi}.$$
Thus, there exists a module idempotent $Q \in \clb(\clr \oplus
\clr)$ with matrix entries in $\clm(\clr)$ such that $Q(\Theta f +
g) =g$ for $f \in \clr$ and $g \in \mbox{ker~}R_{\Psi}$. Moreover,
$\mbox{ran}\, M_{\Theta} = \mbox{ker}\, Q$ and $\mbox{ker}\,
R_{\Psi} = \mbox{ran}\, Q$. The composition $Q \circ
\pi_{\Theta}^{-1} : \clr_{\Theta} \raro \clr$ is well-defined and is
the required invertible module map establishing the similarity of
$\clr_\Theta$ and $\clr$.  \qed

\begin{Corollary}\label{sim2}
Let $\theta_1, \theta_2 \in \clm(H^2_n)$ satisfy $|\theta_1(\z)|^2 +
|\theta_2(\z)|^2 \geq \epsilon$ for all $\z \in \mathbb{B}^n$ and
some $\epsilon >0$. Then the quotient Hilbert module
$(H^2_n)_{\Theta}=(H^2_n \otimes \mathbb{C}^2) / M_{\Theta} H^2_n$
is similar to $H^2_n$.
\end{Corollary}
\NI\textsf{Proof.} The corollary follows from Theorem \ref{sim1}
using the corona theorem for $\clm(H^2_n)$ (see \cite{CSW} or
\cite{OF}).\qed

\subsection{Left invertible multipliers}\label{XX}

The question of similarity of a quotient Hilbert module to the
building block Hilbert module can be raised in the context of a
split short exact sequence. More precisely, suppose $\clr,
\tilde{\clr} \subseteq \clo(\Omega, \mathbb{C})$ be reproducing
kernel Hilbert modules and $\clm(\clr) = \clm(\tilde{\clr})$.
Moreover, suppose $\Theta \in \clm_{\clb(\cle, \cle_*)}(\clr)$ and
hence $\Theta \in \clm_{\clb(\cle, \cle_*)}(\tilde{\clr})$. Consider
the generalized canonical models \[\clr \otimes \cle
\stackrel{M_{\Theta}}\longrightarrow \clr \otimes \cle_*
\longrightarrow \clr_{\Theta} \longrightarrow 0,\quad
\mbox{and}\quad \tilde{\clr} \otimes \cle
\stackrel{M_{\Theta}}\longrightarrow \tilde{\clr} \otimes \cle_*
\longrightarrow \tilde{\clr}_{\Theta} \longrightarrow 0.\]One can
propose the following assertion: $\clr_{\Theta}$ is similar to $\clr
\otimes \clf$ for some Hilbert space $\clf$ if and only if
$\tilde{\clr}_{\Theta}$ is similar to $\tilde{\clr} \otimes
\tilde{\clf}$ for some Hilbert space $\tilde{\clf}$.

Where the answer to the above question is not affirmative, however,
the following hold:

\begin{Theorem}\label{bundle-map}
Let $\clh, \tilde{\clh} \in B_1^*(\Omega)$ for $\Omega \subseteq
\mathbb{C}^n$, be such that $\clm(\clh) \subseteq
\clm(\tilde{\clh})$ and let $\Theta \in \clm_{\clb(\mathbb{C}^p,
\mathbb{C}^q)}(\clh)$ , for $1 \leq p <q$, be left invertible. Then
the similarity of ${\clh}_{\Theta} = (\clh \otimes \mathbb{C}^q)/
M_{\Theta} (\clh \otimes \mathbb{C}^p)$ to $\clh \otimes
\mathbb{C}^{q-p}$ implies the similarity of $\tilde{{\clh}}_{\Theta}
= ({\tilde{\clh}} \otimes \mathbb{C}^q)/ M_{\Theta} ({\tilde{\clh}}
\otimes \mathbb{C}^p)$ to ${\tilde{\clh}} \otimes \mathbb{C}^{q-p}$.
\end{Theorem}

\NI\textsf{Proof.} Since $\clm(\clh) \subseteq \clm(\tilde{\clh})$,
$\Theta \in \clm_{\clb(\mathbb{C}^p, \mathbb{C}^q)}(\clh)$ and
$\tilde{\clh}_{\Theta}$ is well-defined. Moreover, by Theorem
\ref{n-curv}, we have $\clh_{\Theta}, \tilde{\clh}_{\Theta} \in
B^*_{q-p}(\Omega)$. Let $\sigma_{\Theta}$ be a module cross-section
for $\clh_{\Theta}$; that is, $\sigma_{\Theta}: \clh_{\Theta}
\rightarrow \clh \otimes \mathbb{C}^q$ such that
$\pi_{\Theta}\sigma_{\Theta}=I_{\clh_{\Theta}}$. Since $Q :=
\sigma_{\Theta}\pi_{\Theta}$ is a module idempotent on $\clh \otimes
\mathbb{C}^q$, it follows that \[\mbox{ran~}Q \stackrel{.}+
\mbox{ran~}M_{\Theta}= \clh \otimes \mathbb{C}^q.\] But there exists
a $\Phi \in \clm_{\clb(\mathbb{C}^q)}(\clh)$ such that
$$
M_{\Phi}=Q,
$$
and $\Phi(z)$ is an idempotent on $\clb(\mathbb{C}^q)$ for $z \in
\Omega$. An easy argument using localization shows that
$$
\mbox{ran~}\Phi(z) \stackrel{.}+ \mbox{ran~}\Theta(z)= \mathbb{C}^q,
$$
for $z \in \Omega$. But this fact is independent of $\clh$.

\NI Therefore
$$
\tilde{\sigma}_{\Theta}: = M_{\Phi}\tilde{\pi}^{-1}_{\Theta}
$$
is a module map from $\tilde{\clh}_{\Theta}$ to $\tilde{\clh}
\otimes \mathbb{C}^q$, where $\tilde{\pi}_{\Theta}$ is the quotient
map of the short exact sequence for ${\tilde{\clh}}_{\Theta}$.
Moreover, the idempotent
$\tilde{Q}=\tilde{\sigma}_{\Theta}\tilde{\pi}_{\Theta}$ is again
represented by $M_{\Phi}$.

Suppose that $\clh_{\Theta}$ is similar to $\clh \otimes
\mathbb{C}^{q-p}$. Then there exists an invertible module map $X:
\clh \otimes \mathbb{C}^{q-p} \rightarrow \clh_{\Theta}$. Compose
the module maps $\sigma_{\Theta}$ and $X$ to obtain
$Y=\sigma_{\Theta}X: \clh \otimes \mathbb{C}^{q-p} \rightarrow \clh
\otimes \mathbb{C}^q$ and let $\Gamma \in
\clm_{\clb(\mathbb{C}^{q-p}, \mathbb{C}^q)}(\clh)$ so that
$Y=M_{\Gamma}$. Since $\clm(\clh) \subseteq \clm(\tilde{\clh})$, one
can use $\Gamma$ to define
$$
M_{\Gamma}: \tilde{\clh} \otimes \mathbb{C}^{q-p} \rightarrow
\tilde{\clh} \otimes \mathbb{C}^q.
$$
Composing $M^{-1}_{\Gamma}$ and $\tilde{\sigma}_{\Theta}$, one gets
an invertible module map $M^{-1}_{\Gamma} \tilde{\sigma}_{\Theta}$
from $\tilde{\clh}_{\Theta}$ to $\tilde{\clh} \otimes
\mathbb{C}^{q-p}$, which shows that $\tilde{\clh}_{\Theta}$ is
similar to $\tilde{\clh} \otimes \mathbb{C}^{q-p}$. \qed

Theorem \ref{bundle-map} yields the following corollary.
\begin{Corollary}\label{similar}
Let ${\clh} \in B^*_1(\mathbb{D})$ be a contractive Hilbert module
over $A(\mathbb{D})$ and $\Theta \in H^{\infty}_{\clb(\mathbb{C}^p,
\mathbb{C}^q)}(\mathbb{D})$ is left-invertible. Then
${\clh}_{\Theta}$ is similar to ${\clh} \otimes \mathbb{C}^{q-p}$.
\end{Corollary}
\NI\textsf{Proof.} Here one can use Theorem \ref{bundle-map} with
$\tilde{\clh} = H^2(\mathbb{D})$ and ${\clh}$ the given Hilbert
module. Clearly $\clm(\tilde{\clh})=H^{\infty}(\mathbb{D})$. The
proof is completed by appealing to a result of Sz.-Nagy and Foias
about a left invertible $\Theta$ (cf. \cite{NF}). \qed

\NI \textbf{Further results and comments:}

\begin{enumerate}

\item Results of this section can be found in \cite{DKKS2}.

\item The question of similarity is equivalent to a problem in complex
geometry (cf. \cite{D2}). In general, for a split short exact
sequence
$$0 \longrightarrow \clh \otimes \mathbb{C}^p \stackrel{M_\Theta} {\longrightarrow}
\clh \otimes \mathbb{C}^q \stackrel{\pi_{\Theta}} \longrightarrow
  \clh_{\Theta} \longrightarrow 0,$$
one can define the idempotent function $\Gamma: \Omega \rightarrow
\clb(\mathbb{C}^q)$, where $\mbox{ran~}\Gamma$ yields a hermitian
holomorphic subbundle $\clf$ of the trivial bundle $\Omega \times
\mathbb{C}^q$. If $\Omega$ is contractible, then $\clf$ is trivial.
The question of similarity is equivalent to whether one can find a
trivializing frame for which the corresponding Gramian $G$ is
uniformly bounded above and below when
$\clm(\clh)=H^{\infty}(\Omega)$, or  it and its inverse lie in the
multiplier algebra when it is smaller. As mentioned above, this
question is related to the corona problem and the commutant lifting
theorem.

\end{enumerate}

\section{Free resolutions of Hilbert modules}\label{FRHM}

Consideration of dilations such as those in Section 8 in
\cite{JS-HB}, raises the question of what kind of resolutions exist
for co-spherically contractive Hilbert modules over
$\mathbb{C}[\z]$. In particular, Theorem 8.3 in \cite{JS-HB} yields
a unique resolution of an arbitrary pure co-spherically contractive
Hilbert module $\clh$ over $\mathbb{C}[\z]$ in terms of
Drury-Arveson modules and inner multipliers.

Let $\clh$ be a co-spherically contractive Hilbert module over
$\mathbb{C}[\z]$. By Theorem 8.3 in \cite{JS-HB}, there exists a
unique (assuming the minimality) free module $\clf_0$ such that
$\clf_0$ is a dilation of $\clh$. That is, there exists a module
co-isometry \[\pi_{\clh} = X_0 : \clf_0 \raro \clh.\] The kernel of
$X_0$ is a closed submodule of $\clf_0$ and again by Theorem 8.3 in
\cite{JS-HB}, there exists a free module $\clf_1$ and a partially
isometric module map $X_1 : \clf_1 \raro \clf_0$ such that
\[\mbox{ran} X_1 = \mbox{ker~} X_0.\]By repeating this process, one
obtains a sequence of free Hilbert modules $\{\clf_i\}$ and an exact
sequence: \begin{equation}\label{free-res}\cdots
\stackrel{X_2}\longrightarrow \clf_1 \stackrel{X_1} \longrightarrow
\clf_0 \stackrel{X_0} \longrightarrow \clh \longrightarrow
0.,\end{equation} A basic question is whether such a resolution can
have finite length or, equivalently, whether one can take $\cle_N =
\{0\}$ for some finite $N$. That will be the case if and only if
some $X_k$ is an isometry or, equivalently, if $\mbox{ker} X_k =
\{0\}$.

\subsection{Isometric multipliers}

Let $V \in \clb(\cle, \cle_*)$ be an isometry. Then $\Phi_V \in
\clm_{\clb(\cle, \cle_*)}(H^2_n)$, defined by $\Phi_V = I_{H^2_n}
\otimes V$ is an isometric multiplier. The purpose of this
subsection is to prove that all isometric multipliers are of this
form.

\begin{Theorem}\label{res}
For $n > 1$, if $V : H^2_n \otimes \cle \raro H^2_n \otimes \cle_*$
is an isometric module map for Hilbert spaces $\cle$ and $\cle_*$,
then there exists an isometry $V_0 : \cle \raro \cle_*$ such that
$$V( \bm{z}^{\bm{k}} \otimes x)  = \bm{z}^{\bm{k}} \otimes V_0 x,
\;\; \mbox{for~} \bm{k} \in \mathbb{N}^n, x \in \cle_*.$$Moreover,
$\mbox{ran~} V$ is a reducing submodule of $H^2_n \otimes \cle_*$ of
the form $H^2_n \otimes (\mbox{ran}\, V_0)$.
\end{Theorem}
\NI \textsf{Proof.} For $x \in \cle$, $\|x\| = 1$, one can compute
$$V (1 \otimes x) = f(\bm{z}) = \sum_{\bm{k} \in \mathbb{N}^n}
a_{\bm{k}} \bm{z}^{\bm{k}}, \quad \quad \mbox{for} \, \{a_k\}
\subseteq \cle.$$ Then $$V (z_1 \otimes x) = V M_{z_1} (1 \otimes x)
= M_{z_1} V (1 \otimes x) = M_{z_1} f = z_1 f,$$ and $$\|z_1 f \|^2
= \|z_1 V(1 \otimes x)\|^2 = \|z_1 \otimes x\|^2 = 1 =
\|f\|^2.$$Therefore$$\sum_{\bm{k} \in \mathbb{N}^n}
\|a_{\bm{k}}\|^2_{\cle_*} \|\bm{z}^{\bm{k}}\|^2 = \sum_{\bm{k} \in
\mathbb{N}^n} \|a_{\bm{k}}\|^2_{\cle_*} \|\bm{z}^{\bm{k}+e_1}\|^2,
\; \; \mbox{where~} \bm{k} + e_1 = (k_1+1, \ldots, k_n),$$ or
$$\sum_{\bm{k} \in \mathbb{N}^n} \|a_{\bm{k}}\|^2_{\cle_*} \{
\|\bm{z}^{\bm{k} + e_1}\|^2 - \|\bm{z}^{\bm{k}}\|^2\} = 0.$$ If
$\bm{k} = (k_1, \ldots k_n)$, then
\[
\begin{split}
\|\bm{z}^{\bm{k} + e_1}\|^2 & = \frac{(k_1 + 1)! \cdots k_n!}{(k_1 +
\cdots + k_n +1)!} = {\frac{k_1! \cdots k_n!}{(k_1 + \cdots +
k_n)!}} \frac{k_1 + 1}{k_1 + \cdots + k_n + 1}\\ & < {\frac{k_1!
\cdots k_n!}{(k_1 + \cdots + k_n)!}} = \|\bm{z}^{\bm{k}}\|^2,
\end{split}
\]
unless $k_2 = k_3 = \ldots = k_n = 0$. Since, $a_{\bm{k}} \neq 0$
implies $\|\bm{z}^{\bm{k}+e_1}\| = \|\bm{z}^{\bm{k}}\|$ we have $k_2
= \cdots = k_n = 0$. Repeating this argument using $i = 2, \ldots,
n$, it follows that $a_{\bm{k}} = 0$ unless $\bm{k} = (0, \ldots,
0)$ and therefore, $f(\bm{z}) = 1 \otimes y$ for some $y \in
\cle_*$. Set $V_0 x = y$ to complete the first part of the proof.

Finally, since $\mbox{ran~} V = H^2_n \otimes (\mbox{ran} V_0)$, it
follows that $\mbox{ran~} V$ is a reducing submodule, which
completes the proof. \qed

\subsection{Inner resolutions}
This subsection begins with a definition based on the dilation
result in Section 8 in \cite{JS-HB}.

An \textit{inner resolution} of length $N$, for $N = 1, 2, 3,
\ldots, \infty$, for a pure co-spherical contractive Hilbert module
$\clh$ is given by a collection of Hilbert spaces
$\{\cle_k\}_{k=0}^N$, inner multipliers $\varphi_k \in
\clm_{\clb(\cle_k, \cle_{k-1})}(H^2_n)$ for $k = 1, \ldots, N$ with
$X_k = M_{\varphi_k}$ and a co-isometric module map $X_0 : H^2_n
\otimes \cle_0 \raro \clh$ so that
$$\mbox{ran} \, X_k = \mbox{ker}\, X_{k-1},$$ for $k = 0, 1, \ldots,
N$. To be more precise, for $N < \infty$ one has the finite
resolution $$ 0 \longrightarrow H^2_n \otimes \cle_N
\stackrel{X_N}{\longrightarrow} H^2_n \otimes \cle_{N-1}
\longrightarrow \cdots \longrightarrow H^2_n \otimes \cle_1
\stackrel{X_1} \longrightarrow H^2_n \otimes \cle_0 \stackrel{X_0}
\longrightarrow \clh \longrightarrow 0,$$ and for $N = \infty$, the
infinite resolution
$$\cdots \longrightarrow H^2_n \otimes \cle_N \stackrel{X_{N}}{\longrightarrow}
H^2_n \otimes \cle_{N-1} \longrightarrow \cdots \longrightarrow
H^2_n \otimes \cle_1 \stackrel{X_{1}} \longrightarrow H^2_n \otimes
\cle_0 \stackrel{X_{0}} \longrightarrow \clh \longrightarrow 0.$$

The following result shows that an inner resolution does not stop
when $n >1$, unless $\clh$ is a Drury-Arveson module and the
resolution is a trivial one. In particular, the resolution in
(\ref{free-res}) never stops unless $\clh = H^2_n \otimes \cle$ for
some Hilbert space $\cle$.

\begin{Theorem}\label{trivial-res}
If the pure, co-spherically contractive Hilbert module $\clm$
possesses a finite inner resolution, then $\clh$ is isometrically
isomorphic to $H^2_n \otimes \clf$ for some Hilbert space $\clf$.
\end{Theorem}

\NI\textsf{Proof.} By applying Theorem \ref{res} to $X_N$, one can
decompose
\[\cle_{N-1} = \cle^1_{N-1} \oplus \cle^2_{N-1},\] so that
\[\tilde{X}_{N-1} = X_{N-1}|_{H^2_n \otimes \cle^2_{N-1}} \in
\cll(H^2_n \otimes \cle^2_{N-1}, H^2_n \otimes \cle_{N-2}),\] is an
isometry onto $\mbox{ran} \, X_{N-1}$. Hence, one can apply the same
theorem to $\tilde{X}_{N-1}$. Therefore, the desired conclusion
follows by induction. \qed

The following statement proceeds directly from the theorem.

\begin{Corollary}
If $\Theta\in \clm_{\clb(\cle, \cle_*)}(H^2_n)$ is an inner
multiplier for the Hilbert spaces $\cle$ and $\cle_*$ with
$\mbox{ker~}M_{\Theta} = \{0\}$, then the quotient module
$\clh_{\Theta} = (H^2_n \otimes \cle_*)/ \, \mbox{ran} \,M_{\Theta}$
is isometrically isomorphic to $H^2_n \otimes \clf$ for  a Hilbert
space $\clf$. Moreover, $\clf$ can be identified with $(\mbox{ran}\,
V_0)^{\perp}$, where $V_0$ is the isometry from $\cle$ to $\cle_*$
given in Theorem \ref{res}.
\end{Corollary}

Note that in the preceding corollary, one has $\mbox{dim~} \cle_* =
\mbox{dim~} \cle + \mbox{dim~}\clf$.

A resolution of $\clh$ can always be made longer in a trivial way.
Suppose we have the resolution $$ 0 \longrightarrow H^2_n \otimes
\cle_N \stackrel{X_N}{\longrightarrow} H^2_n\otimes \cle_{N-1}
\longrightarrow \cdots \longrightarrow H^2_n \otimes \cle_0
\stackrel{X_0}\longrightarrow \clh \longrightarrow 0.$$ If
$\cle_{N+1}$ is a nontrivial Hilbert space, then define $X_{N+1}$ as
the inclusion map of $H^2_n \otimes \cle_{N+1} \subseteq H^2_n
\otimes (\cle_N \oplus \cle_{N+1})$. Further, set $\tilde{X}_N$
equal to $X_N$ on $H^2_n \otimes \cle_{N} \subseteq H^2_n \otimes
(\cle_{N+1} \oplus \cle_N)$ and equal to $0$ on $H^2_n \otimes
\cle_{N+1} \subseteq H^2_n \otimes (\cle_N \oplus \cle_{N+1})$.
Extending $\tilde{X}_{N}$ to all of $H^2_n \otimes \cle_{N+1}$
linearly, one obtains a longer resolution essentially equivalent to
the original one $$ 0 \longrightarrow H^2_n \otimes \cle_{N+1}
\stackrel{X_{N+1}}{\longrightarrow} H^2_n \otimes (\cle_{N+1} \oplus
\cle_N)  \stackrel{\tilde{X}_N} \longrightarrow \cdots
\longrightarrow  \clh \longrightarrow 0.$$ Moreover, the new
resolution will be inner if the original one is.

The proof of the preceding theorem shows that any finite inner
resolution by Drury-Arveson modules is equivalent to a series of
such trivial extensions of the resolution $$ 0 \longrightarrow H^2_n
\otimes \cle \stackrel{X}{\longrightarrow} H^2_n \otimes \cle
\longrightarrow 0,$$ for some Hilbert space $\cle$ and $X = I_{H^2_n
\otimes \cle}$. Such a resolution will be referred as {\it trivial
inner resolution}. The proof of the following statement is now
straightforward.

\begin{Corollary}
All finite inner resolutions for a pure co-spherically contractive
Hilbert module $\clh$ are trivial inner resolutions.
\end{Corollary}

\subsection{Localizations of free resolutions}
Let $\varphi \in Aut(\mathbb{B}^n)$ and $\varphi = (\varphi_1,
\ldots, \varphi_n)$ where $\varphi_i : \mathbb{B}^n \raro
\mathbb{D}$ is the $i$-th coordinate function of $\varphi$ and $1
\leq i \leq n$. Denote $(H^2_n)_{\varphi}$ by the Hilbert module
\[\mathbb{C}[\z] \times H^2_n \raro H^2_n, \quad \quad (p, f) \mapsto
p(\varphi_1(M_z), \ldots, \varphi_n(M_z))f. \quad \quad (p \in
\mathbb{C}[\z], f \in H^2_n)\]One can check that $(H^2_n)_{\varphi}$
is a co-spherically contractive Hilbert module over
$\mathbb{C}[\z]$. Moreover, as in $n = 1$ case, $H^2_n \cong
(H^2_n)_{\varphi}$ (see \cite{Green}) for all $\varphi \in
Aut(\mathbb{B}^n)$.

In \cite{Green}, D. Green proved the following surprising theorem.

\begin{Theorem}\label{Green-THM}
Let $\clh$ be a co-spherically contractive Hilbert module over
$\mathbb{C}[\z]$ and $\varphi \in Aut(\mathbb{B}^n)$ with $\w =
\varphi^{-1}(0)$ and $\w \in \mathbb{B}^n$. Let (\ref{free-res}) be
the free resolution of $\clh$ with $\clf_i = H^2_n(\cle_i) \oplus
\cls_i$ for some Hilbert space $\cle_i$ and spherical Hilbert module
$\cls_i$ ($i \geq 0$). Then the homology of
\[\cdots \stackrel{X_3(\w)}\longrightarrow \cle_2 \stackrel{X_2(\w)}
\longrightarrow \cle_1 \stackrel{X_1(\w)} \longrightarrow
\cle_0,\]the localization of the free resolution of $\clh$ at $\w
\in \mathbb{B}^n$, is isomorphic to the homology of
\[K((\clh)_{\varphi}): 0 \longrightarrow \cle^n_n((\clh)_{\varphi}) \stackrel{\partial_{n,
(\clh)_{\varphi}}}\longrightarrow \cle^n_{n-1}((\clh)_{\varphi})
\stackrel{\partial_{n-1, (\clh)_{\varphi}}}\longrightarrow \cdots
\stackrel{\partial_{1, (\clh)_{\varphi}}}\longrightarrow
\cle^n_1((\clh)_{\varphi}) \longrightarrow  0,\] the Koszul complex
of $(\clh)_{\varphi}$. Therefore, for all $i \geq 1$ we have
\[\mbox{ker~} \partial_{i, (\clh)_{\varphi}}/ \mbox{ran~} \partial_{i+1,
(\clh)_{\varphi}} \cong \mbox{ker~} X_i(\w)/ \mbox{ran~}
X_{i+1}(\w),\]for each $\w \in \mathbb{B}^n$ and $\varphi \in
Aut(\mathbb{B}^n)$ such that $\varphi(\w) = 0$.
\end{Theorem}

The following result is an immediate consequence of Theorem
\ref{Green-THM}.

\begin{Corollary}\label{Green-c}
Let \[\cdots \stackrel{X_3(\w)}\longrightarrow \cle_2
\stackrel{X_2(\w)} \longrightarrow \cle_1 \stackrel{X_1(\w)}
\longrightarrow \cle_0,\] be the localization at $\w \in
\mathbb{B}^n$ of the free resolution (\ref{free-res}) of a
co-spherically contractive Hilbert module $\clh$ over
$\mathbb{C}[\z]$. Then for all $i \geq n+1$,\[\mbox{ker~} X_{i}(\w)
= \mbox{ran~} X_{i+1}(\w).\]
\end{Corollary}

\NI \textbf{Further results and comments:}

\begin{enumerate}
\item What happens when one relaxes the conditions on the module maps
$\{X_k\}$ so that $\mbox{ran}\, X_{k} = \mbox{ker}\, X_{{k-1}}$ for
all $k$ but do not require them to be partial isometries? In this
case, non-trivial finite resolutions do exist, completely analogous
to what happens for the case of the Hardy or Bergman modules over
$\mathbb{C}[\z]$ for $m>1$. Here is one simple example:

\NI Consider the module $\mathbb{C}_{(0,0)}$ over $\mathbb{C}[z_1,
z_2]$ and the resolution:\[ 0 \longrightarrow H^2_2
\stackrel{X_2}{\longrightarrow} H^2_2 \oplus H^2_2 \stackrel{X_1}
\longrightarrow H^2_2 \stackrel{X_0} \longrightarrow
\mathbb{C}_{(0,0)} \longrightarrow 0,\] where $X_0 f = f(0,0)$ for
$f \in H^2_2$, $X_1(f_1 \oplus f_2) = M_{z_1} f_1 + M_{z_2} f_2$ for
$f_1 \oplus f_2 \in H^2_2 \oplus H^2_2$, and $X_2 f = M_{z_2} f
\oplus ( - M_{z_1} f)$ for $f \in H^2_2$. One can show that this
sequence, which is closely related to the Koszul complex, is exact
and non-trivial; in particular, it does not split as trivial
resolutions do.

\item It is not known if there exists any relationship between the
inner resolution for a pure co-spherically contractive Hilbert
module and more general, {\it not necessarily inner}, resolutions by
Drury-Arveson modules. In particular, is there any relation between
the minimal length of a not necessarily inner resolution and the
inner resolution. Theorem \ref{incl} and Corollary \ref{Green-c}
provides some information on this matter.

\item A parallel notion of resolution for Hilbert modules was studied
by Arveson \cite{A04}, \cite{ATAMS}, which is different from the one
considered in this section. For Arveson, the key issue is the
behavior of the resolution at $0 \in \mathbb{B}^n$ or the
localization of the sequence of connecting maps at $0$. His main
goal, which he accomplishes and is quite non trivial, is to extend
an analogue of the Hilbert's syzygy theorem. In particular, he
exhibits a resolution of Hilbert modules in his class which ends in
finitely many steps.

\item The resolutions considered in (\cite{DM1}, \cite{DM2}) and this
section are related to dilation theory although the requirement that
the connecting maps are partial isometries is sometimes relaxed.

\item Theorem \ref{res} is related to an earlier result of Guo, Hu and
Xu \cite{GHX}.

\item Theorem \ref{Green-THM} and Corollary \ref{Green-c} are due to
Green \cite{Green}. Except that, most of the material is from
\cite{DFS}. However, Theorem \ref{res} was first proved by Arias
\cite{Arias}.
\end{enumerate}

\section{Rigidity}\label{R}

Let $\clh$ be a Hilbert module over $A(\Omega)$ (or, over
$\mathbb{C}[\z]$). Denote by $\clr(\clh)$ the set of all
non-unitarily equivalent submodules of $\clh$, that is, if $\cls_1,
\cls_2 \in \clr(\clh)$ and that $\cls_1 \cong \cls_2$ then $\cls_1 =
\cls_2$.

\NI \textsf{Problem:} Determine $\clr(\clh)$.

By virtue of the characterization results by Beurling and Richter
(see Section 6 in \cite{JS-HB}), we have
\[\clr(H^2(\mathbb{D})) = \{\{0\}, H^2(\mathbb{D})\},  \quad \quad
\mbox{and} \quad \quad \clr(L^2_a(\mathbb{D})) = \{ \cls \subseteq
L^2_a(\mathbb{D}) : \cls \mbox{~is a submodule}\}.\] A Hilbert
module $\clh$ over $A(\Omega)$ is said to be \textit{rigid} if \[
\clr(\clh) = \{ \cls \subseteq \clh : \cls \mbox{~is a submodule}\}
= \mbox{Lat}(\clh).\] Therefore, the Bergman module
$L^2_a(\mathbb{D})$ is rigid. For the Hardy space
$H^2(\mathbb{D}^n)$ with  $n
> 1$
\[\{\{0\}, H^2(\mathbb{D}^n)\} \subset \clr(H^2(\mathbb{D}^n))
\subset \mbox{Lat}(H^2(\mathbb{D}^n))\] The purpose of this section
is to discuss some rigidity results for reproducing kernel Hilbert
modules over $\mathbb{B}^n$ and $\mathbb{D}^n$. For the rest of the
section, unless otherwise stated, it is assumed that $n > 1$.

\subsection{Rigidity of $H^2_n$}
In \cite{GHX}, Guo, Hu and Yu proved that two nested unitarily
equivalent submodules of $H^2_n$ must be equal.
\begin{Theorem}\label{guo-rigid}
Let $\cls_1$ and $\cls_2$ be two submodules of $H^2_n$ and $\cls_1
\subseteq \cls_2$. Then $\cls_1 \cong \cls_2$ if and only if $\cls_1
= \cls_2$.
\end{Theorem}
\NI It is not known whether there exists proper submodules $\cls_1$
and $\cls_2$ of $H^2_n$ such that $\cls_1 \cong \cls_2$ but $\cls_1
\neq \cls_2$. Recall that a submodule $\cls$ of $H^2_n$ is said to
be proper if $\cls \neq H^2_n$, or equivalently, $1 \notin \cls$.

The following result provides a rather weaker version of Theorem
\ref{guo-rigid}.
\begin{Theorem}
If $\cls$ is a submodule of $H^2_n$ which is isometrically
isomorphic to $H^2_n$, then $\cls = H^2_n$.
\end{Theorem}
\NI \textbf{Proof.} The result follows directly from Theorem
\ref{res}.\qed

\subsection{Rigidity of $L^2_a(\mu)$}
The purpose of this subsection is to prove that for a class of
measures $\mu$ on the closure of $\Omega$, two submodules of
$L^2_a(\mu)$ are isometrically isomorphic if and only if they are
equal. The subsection will be concluded by considering when two
submodules of a subnormal Hilbert module $\mathcal{M}$ over
$A(\Omega)$ can be isometrically isomorphic.

Let $\mu$ be the measure on $\overline{\mathbb{D}}$ obtained from
the sum of Lebesgue measure on $\partial\mathbb{D}$ and the unit
mass at 0, then $L^2_a(\mu)$ is not a \v Silov module (see
\cite{DP}). However, it is easy to see that the cyclic submodules
generated by $z$ and $z^2$, respectively, are isometrically
isomorphic but distinct. A quick examination suggests the problem is
that $\mu$ assigns positive measure to the intersection of a zero
variety and $\mathbb{D}$. It turns out that if one excludes that
possibility and $L^2(\nu)$ is not a \v Silov module, then distinct
submodules can not be isometrically isomorphic. The proof takes
several steps.
\begin{Lemma}
Let $\nu$ be a probability measure on $\text{clos } \Omega$ and $f$
and $g$ vectors in $L^2_a(\nu)$ so that the cyclic submodules of
$L^2_a(\nu)$, $[f]$ and $[g]$, generated by $f$ and $g$,
respectively, are isometrically isomorphic  with $f$ mapping to $g$.
Then $|f| = |g|$ a.e. $\nu$.
\end{Lemma}

\NI \textbf{Proof.} If the correspondence $Vf=g$ extends to an
isometric module map, then \[\langle z^{\bm{k}} f, z^{\bm{l}}f
\rangle_{L^2_a(\nu)}\break = \langle z^{\bm{k}}g,
z^{\bm{l}}g\rangle_{L^2_a(\nu)},\] for monomials $z^{\bm{k}}$ and
$z^{\bm{l}}$ in $\mathbb{C}[z]$. This implies that
\[
\int\limits_{\text{clos } \Omega} z^{\bm{k}} \bar{z}^{\bm{l}} |f|^2
d\nu = \int\limits_{\text{clos } \Omega} z^{\bm{k}} \bar{z}^{\bm{l}}
|g|^2 d\nu. \quad \quad ({\bm{k}}, {\bm{l}} \in \mathbb{N}^n)
\]Since the linear span of the set $\{z^{\bm{k}} \bar{z}^{\bm{l}} : {\bm{k}}, {\bm{l}} \in \mathbb{N}^n\}$
forms a self-adjoint algebra which separates the points of clos
$\Omega$, it follows that the two measures $|f|^2~d\nu$ and
$|g|^2~d\nu$ are equal or that $|f| = |g|$ a.e. $\nu$.\qed

The following theorem concerns measures for which point evaluation
on $\Omega$ is bounded.

\begin{Theorem}\label{thm1}
Let $\nu$ be a probability measure on $\text{\rm clos } \Omega$ such
that point evaluation is bounded on $L^2_a(\Omega)$ with closed
support properly containing $\partial\Omega$ but such that $\nu(X) =
0$ for $X$ the intersection of $\text{\rm clos } \Omega$ with a zero
variety. If $\cls_1$ and $\cls_2$ are isometrically isomorphic
submodules of $L^2_a(\nu)$, then $\cls_1 = \cls_2$.
\end{Theorem}

\NI\textbf{Proof.} Let $V$ be an isometric module map from $\cls_1$
onto $\cls_2$. For $0\ne f$ in $\cls_1$, let $g = Vf$. Then by the
previous lemma, it follows that $|f| = |g|$ a.e.\ $\nu$. Since
$\partial\Omega$ is contained in the closed support of $\nu$, it
follows that \[|f(\w)| = |g(\w)|. \quad \quad (\w \in
\partial\Omega)\]Since point evaluation is bounded on
$L^2_a(\Omega), f$ and $g$ are holomorphic on $\Omega$. If \[X =
\{\w \in\Omega : f(\w) = 0\},\] then \[\nu(X) = 0,\]which implies
that $\nu(\Omega\backslash X)>0$. Now \[\sup\limits_{\w
\in\Omega\backslash X} |h(\w)|\le 1,\]where
\[h(\w) : =
\frac{g(\w)}{f(\w)}, \quad \quad (\w \in \Omega \backslash X).\]
Since there is $\w_0$ in the support of $\nu$ in $\Omega\backslash
X$ such that $|h(\w_0)| = 1$, one gets $|h(\w)| \equiv 1$ on
$\Omega\backslash X$. Thus there is a constant $e^{i\theta}$ such
that $f=e^{i\theta}g$ on $\Omega$.

\NI Since this holds for every $f$ in $\cls_1$, by considering
$f_1,f_2$ and $f_1+f_2$, it follows that $Vf=e^{i\theta}f$ for all
$f$ in $\cls_1$ and hence $\cls_1 = \cls_2$. \qed

\vspace{0.3in}

This result contains the results of Richter \cite{Rich}, Putinar
\cite{P}, and Guo--Hu--Xu \cite{GHX} since area measure on
$\mathbb{D}$ or volume measure on $\Omega$ satisfies the hypotheses
of the theorem. However, so do the measures for the weighted Bergman
spaces on $\mathbb{D}$ or weighted volume measure on any domain
$\Omega$.

\begin{Corollary}
If $\Omega$ is a bounded domain in $\mathbb{C}^n$ and $\cls_1$ and
$\cls_2$ are isometrically isomorphic submodules of $L^2_a(\Omega)$,
then $\cls_1 = \cls_2$.
\end{Corollary}

\subsection{Rigidity of $H^2(\mathbb{D}^n)$}
The purpose of this subsection is to discuss the rigidity issue for
a simple class of submodules of $H^2(\mathbb{D}^n)$, namely, the
co-doubly commuting submodules. There is an extensive literature on
rigidity phenomenon for submodules of the Hardy module over
$\mathbb{D}^n$. The reader is referred to the book by Chen and Guo
\cite{C-G}, Chapter 3.

The following rigidity result is due to Agrawal, Clark and Douglas
(Corollary 4 in \cite{ADC}. See also \cite{I}).

\begin{Theorem}
Let $\cls_1$ and $\cls_2$ be two submodules of $H^2(\mathbb{D}^n)$,
both of which contain functions independent of $z_i$ for $i = 1,
\ldots, n$. Then $\cls_1$ and $\cls_2$ are unitarily equivalent if
and only if they are equal.
\end{Theorem}

This yields the following results concerning rigidity of co-doubly
commuting submodules of $H^2(\mathbb{D}^n)$ (see Section
\ref{HMOP}).
\begin{Corollary}\label{Cor1}
Let $\cls_{\Theta} = \sum_{i=1}^n \tilde{\Theta}_i
H^2(\mathbb{D})^n$ and $\cls_{\Phi} = \sum_{i=1}^n \tilde{\Phi}_i
H^2(\mathbb{D})^n$ be a pair of submodules of $H^2(\mathbb{D})^n$,
where $\tilde{\Theta}_i(\bm{z}) = \Theta_i(z_i)$ and
$\tilde{\Phi}_i(\bm{z}) = \Phi_i(z_i)$ for inner functions
$\Theta_i, \Phi_i \in H^{\infty}(\mathbb{D})$ and $\bm{z} \in
\mathbb{D}^n$ and $i = 1, \ldots, n$. Then $\cls_{\Theta}$ and
$\cls_{\Phi}$ are unitarily equivalent if and only if $\cls_{\Theta}
= \cls_{\Phi}$.
\end{Corollary}

\NI\textsf{Proof.} Clearly $\tilde{\Theta}_i \in \cls_{\Theta}$ and
$\tilde{\Phi}_i \in \cls_{\Phi}$ are independent of $\{z_1, \cdots,
z_{i-1}, z_{i+1}, \ldots, z_n\}$ for all $i = 1, \ldots, n$.
Therefore, the submodules $\cls_{\Theta}$ and $\cls_{\Phi}$ contains
functions independent of $z_i$ for all $i = 1, \ldots, n$.
Consequently, if $\cls_\Phi$ and $\cls_\Phi$ are unitarily
equivalent then $\cls_\Theta = \cls_\Phi.$ \qed

\begin{Corollary}\label{Cor2}
Let $\cls_{\Theta} = \sum_{i=1}^n \tilde{\Theta}_i
H^2(\mathbb{D})^n$ be a submodules of $H^2(\mathbb{D})^n$, where
$\tilde{\Theta}_i(\bm{z}) = \Theta_i(z_i)$ for inner functions
$\Theta_i \in H^{\infty}(\mathbb{D})$ for all $i = 1, \ldots, n$ and
$\bm{z} \in \mathbb{D}^n$. Then $\cls_{\Theta}$ and
$H^2(\mathbb{D}^n)$ are not unitarily equivalent.
\end{Corollary}

\NI\textsf{Proof.} The result follows from the previous theorem
along with the observation that $\cls_{\Theta}^{\perp} \neq
\{0\}$.\;\; \qed

\NI \textbf{Further results and comments:}

\begin{enumerate}
\item In \cite{DPSY}, Douglas, Paulsen, Sah
and Yan used algebraic localization techniques to obtain general
rigidity results. In particular, under mild restrictions, they
showed that the submodules obtained from the closure of ideals are
equivalent if and only if the ideals coincide. See also
\cite{D-KY90}, \cite{ADC} for related results.

\item Theorem \ref{thm1} is from \cite{DS1}. In connection with this
section see Richter \cite{Rich}, Putinar \cite{P}, and Guo, Hu and
Xu \cite{GHX}. Corollaries \ref{Cor1} and \ref{Cor2} are from
\cite{JS1}.

\item In \cite{Guo1}, Guo used the notion of
characteristic spaces \cite{Guo2} and obtained a complete
classification submodules of $H^2(\mathbb{D}^n)$ and
$H^2(\mathbb{B}^n)$ generated by polynomials.

\item In \cite{I}, Izuchi proved the
following results: Let $\cls_1$ and $\cls_2$ be submodules of
$H^2(\mathbb{D}^n)$.

\NI (1) If $\mbox{dim}(\cls_1 \ominus \cls_2), \, \mbox{dim}(\cls_2
\ominus \cls_1) < \infty$, then $\cls_1 \cong \cls_2$ if and only if
$\cls_1 = \cls_2$.

\NI (2) Let $\varphi_2$ be an outer function and $\cls_2 =
[\varphi_2]$, the principle submodule generated by $\varphi_2$. If
$\cls_1 \cong \cls_2$ then, $\cls_1 = \Theta \cls_2$ for some inner
function $\Theta \in H^\infty(\mathbb{D}^n)$.

\item For complete reference concerning rigidity for analytic Hilbert
modules, the reader is referred to the book by Chen and Guo
\cite{C-G}.

\end{enumerate}

\section{Essentially normal Hilbert modules}

The purpose of this section is to introduce the notion of
essentially normal Hilbert module, emphasizing a few highlights of
the recent developments in the study of Hilbert modules.

\subsection{Introduction to essential normality}
A Hilbert module $\clh$ over $A$, where $A = A(\Omega)$ or
$\mathbb{C}[\z]$, is said to be \textit{essentially reductive} or
\textit{essentially normal} if the cross-commutators \[[M_i^*, M_j]
= M_i^* M_j - M_j M_i^*,\] are in the ideal of compact operators in
$\clh$ for all $1 \leq i, j \leq n$.

There are many natural examples of essentially normal Hilbert
modules. In particular, $H^2_n$, $L^2_a(\mathbb{B}^n)$ and
$H^2(\mathbb{B}^n)$ are essentially normal. However,
$H^2(\mathbb{D}^n)$ and $L^2_a(\mathbb{D}^n)$ are not essentially
normal whenever $n
> 1$.

In \cite{D06}, Douglas proved the following results: Let $\clh$ be
an essentially normal Hilbert module over $A$ and $\cls$ be a
submodule of $\clh$. Then $\cls$ is essentially normal if and only
if the quotient module $\clq := \clh/ \cls$ is essentially normal.

Another variant of this result concerns a relationship between
essentially normal Hilbert modules and resolutions of Hilbert
modules (see Theorem 2.2 in \cite{D06}):

\begin{Theorem}\label{D-R-EN}
Let $\clh$ be a Hilbert module over $A$ with a resolution of Hilbert
modules \[0 \longrightarrow \clf_1 \stackrel{X}\longrightarrow
\clf_2 \stackrel{\pi}\longrightarrow \clh \longrightarrow 0,\]for
some essentially normal Hilbert modules $\clf_1$ and $\clf_2$. Then
$\clh$ is essentially normal.
\end{Theorem}

The preceding results raise questions about essentially normal
submodules.

\NI\textsf{Problem:} Let $\cls$ be a submodule of $\clh$, where
$\clh = H^2_n$ or $H^2(\mathbb{B}^n)$ or $L^2_a(\mathbb{B}^n)$ and
$n
> 1$. Does it follow that $\cls$ is essentially normal?

This is one of the most active research areas in multivariable
operator theory. For instance, if $\cls$ is a submodule of
$L^2_a(\mathbb{B}^n)$ and generated by a polynomial (by Douglas and
Wang \cite{DW}) or a submodule of $H^2_n$ and generated by a
homogeneous polynomial (by Guo and Wang \cite{GW}), then $\cls$ is
$p$-essentially normal for all $p > n$ (see also \cite{FX},
\cite{E11} and \cite{DRS}).

\subsection{Reductive modules}
This subsection continues the study of unitarily equivalent
submodules of Hilbert modules (see Section 6 in \cite{JS-HB}). In
this context the following problem is of interest.: Let $\clr$ be an
essentially normal quasi-free Hilbert module over $A(\Omega)$ for
which there exists a pure unitarily equivalent submodule. Does it
follow that $\clr$ is subnormal?

Now let $\mathcal{R}$ be a quasi-free Hilbert module over
$A(\Omega)$. Then the Hilbert space tensor product
$\mathcal{R}\otimes H^2(\mathbb{D})$ is a quasi-free Hilbert module
over $A(\Omega\times \mathbb{D})$ which clearly contains the pure
isometrically isomorphic submodule $\mathcal{R}\otimes
H^2_0(\mathbb{D})$. Hence, one can say little without some
additional hypothesis for $\Omega$ or $\mathcal{R}$ or both. Under
the assumption of essential normality on $\mathcal{R}$ the following
holds:

\begin{Theorem}\label{thm5}
Let $\mathcal{R}$ be an essentially normal  Hilbert module over
$A(\Omega)$ and $U$ be an isometric module map $U$  on $\mathcal{R}$
such that \[\bigcap\limits^\infty_{k=0} U^k\mathcal{R} = \{0\}.\]
Then $\mathcal{R}$ is subnormal, that is, there exists a normal
(reductive) Hilbert module $\mathcal{N}$ over $A(\Omega)$ with
$\mathcal{R}$ as a submodule.
\end{Theorem}

\NI \textbf{Proof.} As in the proof of Proposition 6.1 in
\cite{JS-HB}, there exists an isometric isomorphism $\Psi$ from
$\mathcal{R}$ onto $H^2_{\clw}(\mathbb{D})$ with \[\clw =
\mathcal{R}\ominus U\mathcal{R},\]and $\varphi_1,\ldots, \varphi_n$
in $H^\infty_{\mathcal{L}(\clw)}(\mathbb{D})$ such that $\Psi$ is a
$\mathbb{C}[\pmb{z}]$-module map relative to the module structure on
$H^2_{\clw}(\mathbb{D})$ defined so that \[z_j\mapsto T_{\varphi_j}.
\quad \quad(j=1,\ldots, n)\] It remains only to prove that the
$n$-tuple $\{\varphi_1(e^{it}),\ldots, \varphi_n(e^{it})\}$ consists
of commuting normal operators for $e^{it}$-a.e.\ on $\mathbb{T}$.
Then $\mathcal{N}$ is $L^2_{\clw}(\mathbb{T})$ with the module
multiplication defined by $z_i\mapsto L_{\varphi_i}$, where
$L_{\varphi_i}$ denotes pointwise multiplication on
$L^2_{\clw}(\mathbb{T})$. Since the $\{\varphi_j(e^{it})\}^n_{j=1}$
are normal and commute, $L^2_{\clw}(\mathbb{T})$ is a reductive
Hilbert module.

\NI The fact that $\mathcal{R}$ is essentially reductive implies
that each $T_{\varphi_i}$ is essentially normal and hence that the
cross-commutators $[T^*_{\varphi_i},T_{\varphi_j}]$ are  compact for
$1\le i,j\le n$. To finish the proof it suffices to show that
$[T^*_{\varphi_i},T_{\varphi_j}]$ compact implies that
$[L^*_{\varphi_i}, L_{\varphi_j}] = 0$ on $L^2_{\clw}(\mathbb{T})$.

\NI Fix $f$ in $H^2_{\clw}(\mathbb{D})$ and let $N$ be a positive
integer. Next observe that
\begin{equation}\label{tag1}
\lim_{N\to\infty} \|(I-P) L^N_z L^*_{\varphi_i}L_{\varphi_j}f\| = 0,
\end{equation}
and
\begin{equation}\label{tag2}
\lim_{N\to\infty}\|(I-P) L^N_zL^*_{\varphi_i}f\| = 0,
\end{equation}
where $P$ is the projection of $L^2_{\clw}(\mathbb{T})$ onto
$H^2_{\clw}(\mathbb{D})$. Consequently
\begin{align*}
\|[T^*_{\varphi_i}, T_{\varphi_j}] M^N_zf\|
&= \|PL^*_{\varphi_i}PL_{\varphi_j}PL^N_zf - PL_{\varphi_j} PL^*_{\varphi_i}PL^N_zf\|\\
&= \|[L^N_zL^*{\varphi_i} L_{\varphi_j}f - (I-P) L^N_z L^*_{\varphi_i} L_{\varphi_j}f]\\
&\quad - [L_{\varphi_j} L^N_z L^*_{\varphi_i}f - L_{\varphi_j}(I-P)
L^N_z L^*_{\varphi_i}f]\|.
\end{align*}
By (\ref{tag1}) and (\ref{tag2}) one gets
\begin{align*}
\lim_{N\to\infty} \|[T^*_{\varphi_i}, T_{\varphi_j}]L^N_zf\| &=
\lim_{N\to\infty}
\|(L^N_z L^*_{\varphi_i}L_{\varphi_j} - L_{\varphi_j}L^N_z L'_{\varphi_i})f\|\\
&= \lim_{N\to\infty} \|L^N_z[L^*_{\varphi_i}, L_{\varphi_j}]f\| =
\|[L^*_{\varphi_i}, L_{\varphi_j}]f\|.
\end{align*}
Since $[T^*_{\varphi_i}, T_{\varphi_j}]$ is compact and the sequence
$\{e^{iNt}f\}$ converges weakly to 0, it follows that
\[\lim\limits_{N\to\infty}\break \|[T^*_{\varphi_i},T_{\varphi_j}]
e^{iNt}f\| = 0.\] Therefore, \[\|[L^*_{\varphi_i}, L_{\varphi_j}]f\|
= 0.\] Finally, the set of vectors $\{e^{-iNt}f\}\colon \ N\ge 0,
f\in H^2_{\clw}(\mathbb{D})\}$ is norm dense in
$L^2_{\clw}(\mathbb{T})$ and \[\|[L^*_{\varphi_i},L_{\varphi_j}]
e^{-iNt}f\| = \|[L^*_{\varphi_i}, L_{\varphi_j}]f\| = 0.\]
Therefore, $[L^*_{\varphi_i}, L_{\varphi_j}] = 0$, which completes
the proof.\qed

The following result is complementary to Theorem 6.1, \cite{JS-HB}.

\begin{Theorem} \label{reductive}
Let $\clm$ be an  essentially reductive, finite rank, quasi-free
Hilbert module over $A(\mathbb{D})$. Let $U$ be a module isometry
such that \[\cap_{k = 0}^{\infty} U^k \clm = \{0\}.\] Then $\clm$ is
unitarily equivalent to $H^2_{\clf}(\mathbb{D})$ for some Hilbert
space $\clf$ with \[\mbox{dim~} \clf = \mbox{rank~} \clm.\]
\end{Theorem}

\NI \textbf{Proof.} As before (cf. Theorem \ref{thm5}) there is an
isometrical isomorphism, $\Psi\colon \
H^2_{\mathcal{F}}(\mathbb{D})\to \mathcal{M}$ such that $U  = \Psi
T_z\Psi^*$ and there exists $\varphi$ in
$H^\infty_{\mathcal{L}(\mathcal{F})}(\mathbb{D})$ such that $M_z =
\Psi T_\varphi \Psi^*$. Further, since $M_z$  is essentially normal
and $M_z-\omega$ is Fredholm for $\omega$ in $\mathbb{D}$, it
follows that $M_z$ is an essential unitary. Finally, this implies
\[T^*_\varphi T_\varphi-I = T_{\varphi^*\varphi-I},\] is compact and
hence $\varphi^*(e^{it}) \varphi(e^{it}) = I$ a.e.\ or $\varphi$ is
an inner function which completes the proof.\qed

\subsection{Essentially doubly commutativity}

Recall that the Hardy module $H^2(\mathbb{D}^n)$ with $n > 1$ is
doubly commuting but not essentially normal. Therefore, a natural
approach to measure a submodule of the Hardy module
$H^2(\mathbb{D}^n)$ from being small is to consider the cross
commutators $[R_{z_i}^*, R_{z_j}]$ for all $1 \leq i < j \leq n$.

It is difficult in general to characterize the class of essentially
doubly commuting submodules of $H^2(\mathbb{D}^n)$. It is even more
complicated to compute the cross-commutators of submodules of
$H^2(\mathbb{D}^n)$. However, that is not the case for co-doubly
commuting submodules \cite{JS1}:

\begin{Theorem}\label{2-commutator}
Let $\cls = \sum_{i=1}^n \tilde{\Theta}_i H^2(\mathbb{D}^n)$ be a
co-doubly commuting submodule of $H^2(\mathbb{D}^n)$, where
$\tilde{\Theta}_i(\bm{z}) = \Theta_i(z_i)$ for all $\bm{z} \in
\mathbb{D}^n$ and each $\Theta_i \in H^{\infty}(\mathbb{D})$ is
either an inner function or the zero function and $1 \leq i \leq n$.
Then for all $1 \leq i < j \leq n$,
\[[R_{z_i}^* , R_{z_j}] = I_{\clq_{{\Theta}_1}} \otimes \cdots
\otimes \underbrace{P_{\clq_{{\Theta}_i}} M_{z}^*|_{{\Theta}_i
H^2(\mathbb{D})}}_{i^{\rm th}} \otimes \cdots \otimes
\underbrace{P_{{\Theta}_j H^2(\mathbb{D})}
M_{z}|_{\clq_{{\Theta}_j}}}_{j^{\rm th}} \otimes \cdots \otimes
I_{\clq_{{\Theta}_n}},\] and \[\|[R_{z_i}^* , R_{z_j}]\| = (1 -
|{\Theta}_i(0)|^2)^{\frac{1}{2}} (1 -
|{\Theta}_j(0)|^2)^{\frac{1}{2}}.\]
\end{Theorem}

\NI \textsf{Proof.} Let $\cls = \sum_{i=1}^n \tilde{\Theta}_i
H^2(\mathbb{D}^n)$, for some one variable inner functions $\Theta_i
\in H^{\infty}(\mathbb{D})$. Let $\tilde{P}_i$ be the orthogonal
projection in $\cll(\cls)$ defined by
\[\tilde{P}_i = M_{\tilde{\Theta}_i} M_{\tilde{\Theta}_i}^*,\] for
all $i = 1, \ldots, n$. By virtue of Corollary \ref{S-proj} and
Lemma \ref{P-F},
\[\begin{split} {P}_{\cls} & = I_{H^2(\mathbb{D}^n)} - \mathop{\Pi}_{i=1}^n
(I_{H^2(\mathbb{D}^n)} - \tilde{P}_i) \\ & = \tilde{P}_1 (I -
\tilde{P}_2) \cdots (I - \tilde{P}_n) + \tilde{P}_2 (I -
\tilde{P}_3) \cdots (I - \tilde{P}_n) + \cdots + \tilde{P}_{n-1} (I
- \tilde{P}_n) + \tilde{P}_n\\ & = \tilde{P}_n (I - \tilde{P}_{n-1})
\cdots (I - \tilde{P}_1) + \tilde{P}_{n-1} (I - \tilde{P}_{n-2})
\cdots (I - \tilde{P}_1) + \cdots + \tilde{P}_2 (I - \tilde{P}_1) +
\tilde{P}_1,
\end{split}\]and
\[{P}_{\clq} = \mathop{\Pi}_{i=1}^n (I_{H^2(\mathbb{D}^n)} - \tilde{P}_i).\] On the other hand, for all $1 \leq i < j \leq n$,
one gets
\[\begin{split}[R_{z_i}^*, R_{z_j}] & = {P}_{\cls} M_{z_i}^* M_{z_j}|_{\cls} -
{P}_{\cls} M_{z_j} {P}_{\cls} M_{z_i}^*|_{\cls},\end{split}
\]and that
\[\begin{split}{P}_{\cls} M_{z_i}^*
M_{z_j}{P}_{\cls} - {P}_{\cls} M_{z_j} {P}_{\cls}
M_{z_i}^*{P}_{\cls} & = {P}_{\cls} M_{z_i}^* M_{z_j}{P}_{\cls} -
{P}_{\cls} M_{z_j} (I - {P}_{\clq}) M_{z_i}^*{P}_{\cls}\\ & =
{P}_{\cls} M_{z_j} {P}_{\clq} M_{z_i}^* {P}_{\cls}.
\end{split} \]Furthermore, for all $1 \leq i < j \leq n$,  \[\begin{split} {P}_{\cls} & M_{z_j} {P}_{\clq} M_{z_i}^*
{P}_{\cls} \\ & = [\tilde{P}_n (I - \tilde{P}_{n-1}) \cdots (I -
\tilde{P}_1) + \tilde{P}_{n-1} (I - \tilde{P}_{n-2}) \cdots (I -
\tilde{P}_1) + \cdots + \tilde{P}_2 (I - \tilde{P}_1) +
\tilde{P}_1]\\& \;\; \;\; M_{z_j} [\mathop{\Pi}_{l=1}^n
(I_{H^2(\mathbb{D}^n)} - \tilde{P}_l)] M_{z_i}^*
\\& \;\;\;\; [\tilde{P}_1 (I - \tilde{P}_2) \cdots (I - \tilde{P}_n) + \tilde{P}_2 (I - \tilde{P}_3) \cdots
(I - \tilde{P}_n) + \cdots + \tilde{P}_{n-1} (I - \tilde{P}_n) +
\tilde{P}_n]\\ & = [\tilde{P}_n (I - \tilde{P}_{n-1}) \cdots (I -
\tilde{P}_1) + \tilde{P}_{n-1} (I - \tilde{P}_{n-2}) \cdots (I -
\tilde{P}_1) + \cdots + \tilde{P}_2 (I - \tilde{P}_1) +
\tilde{P}_1]\\& \;\; \;\; [\mathop{\Pi}_{l \neq j}
(I_{H^2(\mathbb{D}^n)} - \tilde{P}_l)] M_{z_j} M_{z_i}^*
[\mathop{\Pi}_{l\neq i} (I_{H^2(\mathbb{D}^n)} - \tilde{P}_l)] \\&
\;\;\;\; [\tilde{P}_1 (I - \tilde{P}_2) \cdots (I - \tilde{P}_n) +
\tilde{P}_2 (I - \tilde{P}_3) \cdots (I - \tilde{P}_n) + \cdots +
\tilde{P}_{n-1} (I - \tilde{P}_n) + \tilde{P}_n]
\\ & = [\tilde{P}_j (I - \tilde{P}_{j-1}) \cdots
(I - \tilde{P}_1)] M_{z_i}^* M_{z_j} [\tilde{P}_i (I - \tilde{P}_{i+1}) \cdots (I - \tilde{P}_n)]\\
& = [(I - \tilde{P}_1) \cdots (I - \tilde{P}_{j-1}) \tilde{P}_j]
M_{z_i}^* M_{z_j} [\tilde{P}_i (I - \tilde{P}_{i+1}) \cdots (I -
\tilde{P}_n)].
\end{split}\]
These equalities shows that \[\begin{split}[R_{z_i}^*, R_{z_j}] & =
[(I - \tilde{P}_1) \cdots (I - \tilde{P}_i) \cdots (I -
\tilde{P}_{j-1}) \tilde{P}_j] M_{z_i}^* M_{z_j} [\tilde{P}_i (I -
\tilde{P}_{i+1}) \cdots (I - \tilde{P}_j) \cdots (I -
\tilde{P}_n)]\\ & = (I - \tilde{P}_1) (I - \tilde{P}_2) \cdots
(I - \tilde{P}_{i-1}) \;((I - \tilde{P}_i) M_{z_i}^* \tilde{P}_i)\; (I - \tilde{P}_{i+1}) \cdots \\
& \;\;\;\; \cdots(I - \tilde{P}_{j-1}) \;(\tilde{P}_j M_{z_j} (I -
\tilde{P}_j))\; (I - \tilde{P}_{j+1}) \cdots
(I-\tilde{P}_n).\end{split}\] Moreover,
\[[R_{z_i}^*, R_{z_j}] = [(I - \tilde{P}_1) \cdots (I - \tilde{P}_{j-1}) \tilde{P}_j]
M_{z_i}^* M_{z_j} [(I - \tilde{P}_1) \cdots (I - \tilde{P}_{i-1})
\tilde{P}_i (I - \tilde{P}_{i+1}) \cdots (I - \tilde{P}_n)],\] and
\[[R_{z_i}^*, R_{z_j}] = [(I - \tilde{P}_1) \cdots (I - \tilde{P}_{j-1}) \tilde{P}_j (I -
\tilde{P}_{j+1}) \cdots (I - \tilde{P}_n)] M_{z_i}^* M_{z_j}
[\tilde{P}_i (I - \tilde{P}_{i+1}) \cdots (I - \tilde{P}_n)].\] Now
we can conclude that\[\begin{split}[R_{z_i}^*& , R_{z_j}] =
I_{\clq_{\Theta_1}} \otimes \cdots \otimes
\underbrace{{P}_{\clq_{\Theta_i}} M_{z}^*|_{\Theta_i
H^2(\mathbb{D})}}_{i^{\rm th}} \otimes \cdots \otimes
\underbrace{{P}_{\Theta_j H^2(\mathbb{D})}
M_{z}|_{\clq_{\Theta_j}}}_{j^{\rm th}} \otimes \cdots \otimes
I_{\clq_{\Theta_n}}. \end{split}\] Further, note that  \[
\begin{split} \|[R_{z_i}^* , R_{z_j}]\| & = \|I_{\clq_{\Theta_1}}
\otimes \cdots \otimes {P}_{\clq_{\Theta_i}} M_{z}^*|_{\Theta_i
H^2(\mathbb{D})} \otimes \cdots \otimes {P}_{\Theta_j
H^2(\mathbb{D})} M_{z}|_{\clq_{\Theta_j}} \otimes \cdots \otimes
I_{\clq_{\Theta_n}}\| \\& = \|{P}_{\clq_{\Theta_i}}
M_{z}^*|_{\Theta_i H^2(\mathbb{D})}\| \|{P}_{\Theta_j
H^2(\mathbb{D})} M_{z}|_{\clq_{\Theta_j}}\|, \end{split}\] and
consequently by Proposition 2.3 in \cite{JS1} it follows that \[
\|[R_{z_i}^* , R_{z_j}]\| = (1 - |\Theta_i(0)|^2)^{\frac{1}{2}} (1 -
|\Theta_j(0)|^2)^{\frac{1}{2}}.\]This completes the proof. \qed

The following corollary reveals the significance of the identity
operators in the cross commutators of the co-doubly commuting
submodules of $H^2(\mathbb{D}^n)$ for $n >2$.

\begin{Corollary}\label{EN-2}
Let $\cls = \sum_{i=1}^n \tilde{\Theta}_i H^2(\mathbb{D}^n)$ be a
submodule of $H^2(\mathbb{D}^n)$ for some one variable inner
functions $\{\tilde{\Theta}_i\}_{i=1}^n \subseteq
H^{\infty}(\mathbb{D}^n)$. Then

\NI (1) for $n = 2$: the rank of the cross commutator of $\cls$ is
at most one and the Hilbert-Schmidth norm of the cross commutator is
given by
\[\|[R_{z_1}^* , R_{z_2}]\|_{\bm\, HS} = (1 -
|{\Theta}_1(0)|^2)^{\frac{1}{2}} (1 -
|{\Theta}_2(0)|^2)^{\frac{1}{2}}.\]In particular, $\cls$ is
essentially doubly commuting.

\NI (2) for $n >2$: $\cls$ is essentially doubly commuting (or of
Hilbert-Schmidth cross-commutators) if and only if that $\cls$ is of
finite co-dimension, that is,
\[\mbox{dim~} [H^2(\mathbb{D}^n)/\cls] < \infty.\]Moreover, in this case, for all $1 \leq i <
j \leq n$ \[ \|[R_{z_i}^* , R_{z_j}]\|_{\bm\, HS} = (1 -
|\Theta_i(0)|^2)^{\frac{1}{2}} (1 -
|\Theta_j(0)|^2)^{\frac{1}{2}}.\]
\end{Corollary}

The following statements also proceeds directly from the theorem.

\begin{Corollary}\label{<3}
Let $n > 2$ and $\cls = \sum_{i=1}^k \tilde{\Theta}_i
H^2(\mathbb{D}^n)$ be a co-doubly commuting proper submodule of
$H^2(\mathbb{D}^n)$ for some inner functions $\{\Theta_i\}_{i=1}^k$
and $k < n$. Then $\cls$ is not essentially doubly commuting.
\end{Corollary}

\begin{Corollary} Let
$\cls$ be a co-doubly commuting submodule of $H^2(\mathbb{D}^n)$ and
$\clq: = H^2(\mathbb{D}^n)/\cls$ and $n > 2$. Then the following are
equivalent:

(i)  $\cls$ is essentially doubly commuting.

(ii) $\cls$ is of finite co-dimension.

(iii) $\clq$ is essentially normal.
\end{Corollary}

The following one is a "rigidity" type result.

\begin{Corollary}\label{n=2}
Let $n  \geq 2$ and $\cls = \sum_{i=1}^n \tilde{\Theta}_i
H^2(\mathbb{D}^n)$ be an essentially normal co-doubly commuting
submodule of $H^2(\mathbb{D}^n)$ for some one variable inner
functions $\{\Theta_i\}_{i=1}^n$. If $\cls$ is of infinite
co-dimension, then $n = 2$.
\end{Corollary}
\NI\textsf{Proof.}  The result follows from Theorem
\ref{2-commutator} and the fact that a finite co-dimensional
submodule of an essentially doubly commuting Hilbert module over
$\mathbb{C}[\bm{z}]$ is essentially doubly commuting. \qed

It is now clear that the general picture of essentially doubly
commuting submodules of $H^2(\mathbb{D})^n$ is much more
complicated.

\NI \textbf{Further results and comments:}

\begin{enumerate}

\item It is an extremely interesting question as to whether essential
reductivity is related to a lack of corners or not being a product.

\item In \cite{AC}, Ahern and Clark proved that there exists a bijective
correspondence between submodules of $H^2(\mathbb{D}^n)$ of finite
codimension, and the ideals in $\mathbb{C}[\z]$ of finite
codimension whose zero sets are contained in $\mathbb{D}^n$. In
\cite{GZ02}, Guo and Zheng characterized the finite co-dimensional
submodules of the Bergman module and the Hardy module over
$\mathbb{B}^n$ or $\mathbb{D}^n$ (also see Corollary 2.5.4 in
\cite{C-G}).

\begin{Theorem}
Let $\Omega = \mathbb{B}^n$ or $\mathbb{D}^n$ and $\cls$ be a
submodule of $L^2_a(\Omega)$ or $H^2(\Omega)$. Then $\cls$ is of
finite co-dimension if and only if $\cls^\perp$ consists of rational
functions.
\end{Theorem}

\item Second subsection is from \cite{DS1} and the final subsection is from \cite{JS1}.
Part (1) of the Corollary \ref{EN-2} was obtained by R. Yang
(Corollary 1.1, \cite{Y-JOT05}).

\item In \cite{Berger}, Berger and Shaw proved a surprising result
concerning essentially normal Hilbert modules. Suppose $\clh$ be a
\textit{hyponormal} Hilbert module over $\mathbb{C}[z]$, that is,
$[M^*, M] \geq 0$. Moreover, assume that $\clh$ is
\textit{rationally finitely generated}, that is, there exists $m \in
\mathbb{N}$ and $\{f_1, \ldots, f_m\} \subseteq \clh$ such that
\[\{\mathop{\sum}_{i=1}^m r_i(M) f_i : r_i \in Rat(\sigma(M))\}\]is
dense in $\clh$. Then
\[trace[M^*, M] \leq \frac{m}{\pi} \mbox{Area}(\sigma(M)).\]In
particular, every rationally finitely generated hyponormal Hilbert
module is essentially normal.  It is not known whether the
Berger-Shaw theorem holds for "hyponormal" Hilbert modules over
$\mathbb{C}[\z]$. However, in \cite{D-KY}, Douglas and Yan proposed
a version of Berger-Shaw theorem in several variables under the
assumption that the spectrum of the Hilbert module is contained in
an algebraic curve (see also \cite{ZIEOT}). The reader is also
referred to the work of Chavan \cite{Chavan} for a different
approach to the Berger-Shaw theorem in the context of
$2$-hyperexpansive operators.

\item In connection with trace formulae, integral operators,
fundamental trace forms and pseudo-differential operators see also
Pincus \cite{Pin}, Helton and Howe \cite{HeHo} and Carey and Pincus
\cite{CP1}, \cite{CP2}. See also the recent article by Howe
\cite{Howe}.

\item Let $\cls$ be a homogeneous submodule of $H^2(\mathbb{D}^2)$. In \cite{CMY}, Curto, Muhly and Yan proved
that $\cls$ is always essentially doubly commuting.

\item The reader is referred to the work by Ahern and Clark \cite{AC} for more details on finite
co-dimensional submodules of the Hardy modules over $\mathbb{D}^n$
(see also \cite{C-G}).

\item In \cite{AD-LAAA}, Alpay and Dubi characterized finite
co-dimensional subspaces of $H^2_n \otimes \mathbb{C}^m$ for $m \in
\mathbb{N}$ (see also \cite{AD}).

\end{enumerate}

\end{document}